\documentclass[reqno,a4paper,12pt]{amsart}

\usepackage[a4paper,hmargin=2cm,tmargin=3cm,bmargin=3cm]{geometry}
\usepackage[utf8x]{inputenc}
\usepackage[T1]{fontenc}
\usepackage{amsmath,amssymb,amstext,amsthm,amscd,mathrsfs,bm,slashed}
\usepackage{graphicx,color}
\usepackage{array}
\usepackage{cite}
\usepackage{hyperref}
\usepackage{amsfonts}
\usepackage{amscd}
\usepackage{verbatim} 
\usepackage{fancyhdr}
\usepackage{tikz}
\usepackage{tikz-cd}
\usepackage{slashed}
\usetikzlibrary{matrix}
\usetikzlibrary{positioning}
\usepackage[all]{xy}

\usepackage{epsfig}
\usepackage{amssymb}
\usepackage{tabularx}
\usepackage{calligra}
\usepackage{enumerate}
\usepackage{enumitem}
\usepackage{hyperref}
\usepackage{float}
\usepackage{stackrel}
\usepackage{longtable}
\usepackage{mathrsfs}
\usepackage{verbatim} 
\usepackage{units}
\usepackage{textcomp}
\usepackage{turnstile}
\usepackage{ stmaryrd }
\usepackage{tikz-3dplot}

\usepackage[font=small]{caption}

\numberwithin{equation}{section}

\theoremstyle{plain}
\newtheorem{proposition}{Proposition}
\newtheorem{theorem}[proposition]{Theorem}

\def\d{{\rm d}}
\def\i{{\rm i}}
\def\CP{\mathbb{CP}}
\newcommand{\X}{\mathcal{X}}
\newcommand{\Y}{\mathcal{Y}}

\begin{document}

\title{On Chen-Teo geometries with cosmological constant}

\author[B. Araneda]{Bernardo Araneda}
\address{
School of Mathematics and Maxwell Institute for Mathematical Sciences,
University of Edinburgh, EH9 3FD, United Kingdom, 
\newline
Max-Planck-Institut f\"ur Gravitationsphysik (Albert-Einstein-Institut), 
Am M\"uhlenberg 1, D-14476 Potsdam, Germany}
\email{baraneda@ed.ac.uk, bernardo.araneda@aei.mpg.de}
\date{\today}

\maketitle

\begin{abstract}
The Chen-Teo geometry is a Riemannian, Ricci-flat ALF 4-manifold, containing an AF gravitational instanton that gives the first counterexample to the Euclidean black hole uniqueness conjecture. We investigate the problem of constructing an Einstein analogue with a non-zero cosmological constant $\lambda$. We show that the solution is either the Pleba\'nski-Demia\'nski metric with $\lambda$, or it has an anti-self-dual Weyl tensor. We study the latter case in detail: we prove that for $\lambda<0$, there is a conformal infinity separating two asymptotically hyperbolic metrics; we show that one of them is globally conformal to an ALE scalar-flat K\"ahler metric; we construct gravitational instantons with different topologies; and we show that the geometry is a 4-pole solution in the Calderbank-Pedersen classification.
\end{abstract}

\section{Introduction}
\label{sec:introduction}

In this paper, a gravitational instanton refers to a four-dimensional, complete, oriented Riemannian Einstein manifold, with a possibly non-zero cosmological constant $\lambda$. 
If $\lambda=0$, the manifold is Ricci-flat. In this case one is often interested in non-compact instantons with a prescribed global structure, such as asymptotically flat (AF), asymptotically locally flat (ALF), or asymptotically locally Euclidean (ALE), although there are also other possibilities, for more details see e.g. \cite{LiSun2}. 
Standard ALE and ALF examples are the (multi-)Eguchi-Hanson metrics for the former, and the (multi-)Taub-NUT and Taub-bolt instantons for the latter. In the AF case, the classical non-trivial example is the Euclidean Kerr solution (with Euclidean Schwarzschild as a particular case).

The Euclidean Black Hole Uniqueness Conjecture postulated that Kerr is the only non-trivial AF, Ricci-flat instanton, but this was disproved with the discovery of the Chen-Teo instanton \cite{CT1}. 
Different kinds of generalisations of this solution have now been obtained, such as a larger family of (non-smooth) ALF metrics \cite{CT2}, Einstein-Maxwell constructions \cite{a23, AD25}, and an infinite family of new AF Ricci-flat instantons, which are not explicitly known but whose existence was proved by Li and Sun \cite{LiSun1}. In a recent remarkable paper \cite{Teo}, Teo found explicitly the next instanton in the Li-Sun family.
The inclusion of a cosmological constant in the Chen-Teo metric, however, remains an open problem, which motivates this paper.

With the exception of the new Teo instanton \cite{Teo}, all of the explicitly known Ricci-flat metrics mentioned above are Hermitian (i.e., they have a compatible integrable complex structure), and furthermore conformally K\"ahler. If the conformal factor is constant, the metrics can be constructed with the Gibbons-Hawking ansatz. If it is non-constant, the problem is more difficult, but works from Przanowski {\it et al} \cite{P1, P2} and Tod \cite{Tod:2020ual} reduce the Ricci-flat condition to a single scalar, non-linear, integrable PDE, called the $SU(\infty)$ Toda equation.

For a non-zero cosmological constant $\lambda$, the problem can still be reduced to a `modified' Toda equation, but this is no longer an integrable equation  \cite{Tod:2020ual}. This means, in particular, that the existence of a possible `Chen-Teo metric with $\lambda$' is {\it a priori} not guaranteed. However, an argument in favour of the existence of some version of it, at least for $\lambda>0$, can be deduced from the work of Biquard-Gauduchon and LeBrun, as follows. 

For $\lambda>0$, Hermitian-Einstein (non-K\"ahler) instantons were classified by LeBrun \cite{LeBrun2}: these are the (anti) Fubini-Study metric on $\CP^2$, the Page metric on $\CP^2\#\overline{\CP}{}^{2}$, and the Chen-LeBrun-Weber (CLW) metric on $\CP^2\#2\overline{\CP}{}^{2}$. These metrics are toric and conformally K\"ahler. Toric K\"ahler structures can be encoded in Delzant polytopes, which represent the degeneration of the torus action, cf. \cite{Abreu}. 
For $\lambda=0$, Hermitian (non-K\"ahler), ALF, Ricci-flat instantons were classified by Biquard-Gauduchon \cite{BG} and Li \cite{Li23}: the only solutions are anti-Taub-NUT on $\mathbb{R}^4$, Taub-bolt on $\CP^2\setminus\{\rm pt\}$, Kerr on $\mathbb{R}^2\times S^2$, and Chen-Teo on $\CP^2\setminus S^1$. The solutions are again toric and conformally K\"ahler. Biquard-Gauduchon noted an interesting duality between the $\lambda=0$ and $\lambda>0$ cases: the Delzant polytopes of the ALF instantons correspond to those of the compact instantons after removing an appropriate edge of the polytope \cite{BG}; cf. also \cite{OSD}.
In particular, this suggests to interpret the CLW metric on $\CP^2\#2\overline{\CP}{}^{2}$ as a version of the Chen-Teo instanton with $\lambda>0$. The issue, however, is that although the CLW metric is known to exist \cite{CLW}, its explicit form remains unknown.

Thus, despite the non-integrability of the modified Toda equation, the above argument motivates trying to construct the `Chen-Teo metric with $\lambda$'; that is, we consider the ansatz \cite{CT2}
\begin{align}
&g = \frac{(F\d\tau+G\d\phi)^2}{(x-y)FH}+\frac{kH}{(x-y)^3}\left(\frac{\d{x}^2}{\mathcal{X}} - \frac{\d{y}^2}{\mathcal{Y}} - \frac{\mathcal{X}\mathcal{Y}}{kF}\d\phi^2\right), \label{CTansatz} \\
&\X=\X(x), \quad \Y=\Y(y), \quad F=F(x,y), \quad G=G(x,y), \quad H=H(x,y) \label{CTansatzfunctions}
\end{align}
(or a more general but structurally close ansatz, cf. \S\ref{sec:ansatz}) where $k$ is a constant and the functions \eqref{CTansatzfunctions} are {\em arbitrary}, and try to determine these functions so as to obtain a solution to the Einstein equations with $\lambda\neq0$, so that the solution reduces to the Ricci-flat Chen-Teo metric when $\lambda\to0$. In this paper we wish to address this problem. We will, however, find a strong constraint on the possible solutions, summarised as follows:

\begin{theorem}\label{thm:solutions}
If $\lambda\neq0$ and the Weyl tensor is not (anti-)self-dual, the solution is locally isometric to the Pleba\'nski-Demia\'nksi metric with $\lambda$.
\end{theorem}

Therefore, within the natural ansatz considered here, the Einstein equations exhibit a remarkable rigidity: the only non-(anti-)self-dual deformation with $\lambda$ is the Pleba\'nski–Demia\'nski family. In particular, this means that the CLW metric cannot be reached with \eqref{CTansatz}-\eqref{CTansatzfunctions}. In \S\ref{sec:discussion} we will argue that this can actually be expected even before solving the equations, based on an analysis of Delzant polytopes.

The Pleba\'nski-Demia\'nski metric with $\lambda$ is well-known \cite{PD}. But if we instead solve the equations allowing an anti-self-dual (ASD) Weyl tensor, the resulting Einstein metric we find has a number of interesting features that appropriately generalise those of the Ricci-flat Chen-Teo metric to the ASD Einstein case with $\lambda\neq0$ (cf. \cite{CT2}). 

\begin{theorem}\label{thm:ASDCT}
Consider the following four-dimensional Riemannian metric:
\begin{equation}\label{ASDCTintro}
\begin{aligned}
g ={}& \frac{(F\d\tau+G\d\phi)^2}{(x+\nu y)FH}
+\frac{kH}{(x+\nu y)^3}\left(\frac{\d{x}^2}{\X} - \frac{\d y^2}{\Y} - \frac{\X\Y}{kF}\d\phi^2\right), \\
\X ={}& a_0 + a_1 x + a_2 x^2 + a_3 x^3 + a_4 x^4, \quad 
\Y = a_0 + a_1 y + a_2 y^2 + a_3 y^3 + a_4 y^4, \\
F ={}& y^2\X - x^2\Y, \\
H ={}& \frac{3(\nu+1)}{2k\lambda}(x+\nu{y})\left[ 2(\nu-1)(a_0 - a_4 x^2{y}^2) - (x-\nu{y})(a_1 - a_3 x {y})\right], \\
G ={}& \frac{3(\nu+1)}{2k \lambda \nu}\left[ 
\left(-a_3 \nu{y}^3-2(\nu-1)a_2 {y}^2 - 2a_0(\nu-1)-a_1(2\nu-1){y} \right) \X \right. \\
& \qquad\qquad \left. + \left(2a_2(\nu-1)x^2+2a_4(\nu-1)x^4+a_1\nu x+a_3(2\nu-1)x^3 \right)\Y \right],
\end{aligned}
\end{equation}
where $a_0,...,a_4,k,\nu$  and $\lambda<0$ are real constants. Then:
\begin{enumerate}
\item The metric $g$ is Einstein with cosmological constant $\lambda$, and the Weyl tensor is ASD.
\item The solution is asymptotically locally hyperbolic. There is a conformal infinity located at $x+\nu y=0$, with the topology of a lens space.
\label{item:ALH}
\item There is another ASD Einstein metric beyond conformal infinity, which is locally isometric but with a generally different global topology.
\label{item:ASDE2}
\item The metric $g_{K}=\frac{(x+\nu y)^2}{(x-y)^2}g$ is K\"ahler, scalar-flat, and 
ALE of order $2$.
\label{item:SFK}
\item If $\nu\neq\pm1$, the metric \eqref{ASDCTintro} is a 4-pole solution in the Calderbank-Pedersen classification. If $\nu=1$, it reduces to a 3-pole solution, and it is isometric to the ASD Pleba\'nski-Demia\'nski metric with $\lambda$. In the limit $\nu\to-1$, the metric reduces to a 3-centred ALE Gibbons-Hawking space.
\label{item:CP}
\end{enumerate}
\end{theorem}

The Calderbank-Pedersen classification, mentioned in point \eqref{item:CP}, is a local description of all (anti-)self-dual Einstein metrics with toric symmetry, in terms of an eigenfunction $\mathcal{F}$ of the Laplacian in the hyperbolic plane \cite{CP}. A special role in this construction is played by the so-called $m$-pole solutions: those in which $\mathcal{F}$ is of the form
\begin{align}\label{mpole}
\mathcal{F} = \sum_{k=1}^{m}\frac{\sqrt{\alpha_{k}^2\rho^2+(\alpha_k\zeta-\beta_k)^2}}{\sqrt{\rho}},
\end{align}
for some constants $\alpha_k,\beta_k$, where $\rho,\zeta$ are special coordinates adapted to the toric symmetry. We will show that the metric \eqref{ASDCTintro} is a 4-pole solution in this scheme. Since an $m$-pole solution depends on $2m-4$ parameters \cite{CP}, this implies that \eqref{ASDCTintro} depends effectively on 4 parameters\footnote{In particular, one can see from \eqref{ASDCTintro} that $k$ is not a genuine parameter. In our construction of \eqref{ASDCTintro}, $k$ is artificially introduced in \S\ref{sec:ansatz} below.}.
Given the two limits $\nu\to\pm1$, the full solution \eqref{ASDCTintro} can be interpreted as an ASD Einstein metric interpolating between 3-centred ALE Gibbons-Hawking and ASD Pleba\'nski-Demia\'nski with $\lambda$, similarly to the Ricci-flat Chen-Teo metric \cite{CT2}.

Generically, the solution \eqref{ASDCTintro} has conical and orbifold singularities, so it is not a gravitational instanton. To obtain a smooth solution these singularities must be removed, which, in general, may or may not be possible. We describe the situation in terms of the diagrams in Fig. \ref{fig:diagrams}, which represent the domain of the $x,y$ coordinates. They encode the degeneration of the torus action (solid edges and vertices) and the location of conformal infinities (dashed lines). For a given diagram, the two ASD Einstein metrics are denoted `${\rm L}$' (Left) and `${\rm R}$' (Right), in relation to their location with respect to conformal infinity. It turns out that in all cases, the solution `R' is globally conformal to the ALE K\"ahler metric of Theorem \ref{thm:ASDCT}, while the solution `L' is not. We will show that \eqref{ASDCTintro} indeed contains gravitational instantons:

\begin{theorem}\label{thm:instantons}
The solution `R' in Fig. \ref{fig:diagrams} yields a gravitational instanton with:
$(a)$ four continuous parameters, and the topology of the 4-ball $B^4$;
$(b)$ two continuous parameters, and the topology of the total space of the line bundle $\mathcal{O}(-p)\to\mathbb{CP}^1$ for $p\neq0$;
$(c)$ two discrete parameters, and the topology of the toric resolution of $\mathbb{C}^2/\Gamma$ where the exceptional divisor consists of two 2-spheres.
\end{theorem}

 The cyclic subgroup $\Gamma\subset U(2)$ is as usual generated by ${\rm diag}(e^{2\pi i/p},e^{2\pi i q/p})$. The topology has been determined using \cite{CS} and the fact that the associated K\"ahler metric is ALE.

Finally, we mention that classification results for toric ASD Einstein gravitational instantons have been recently given in \cite{ALR26}, building upon \cite{CS}. The gravitational instanton associated to case $(c)$ above corresponds to the three-fixed-points Calderbank-Singer instanton.

\begin{center}
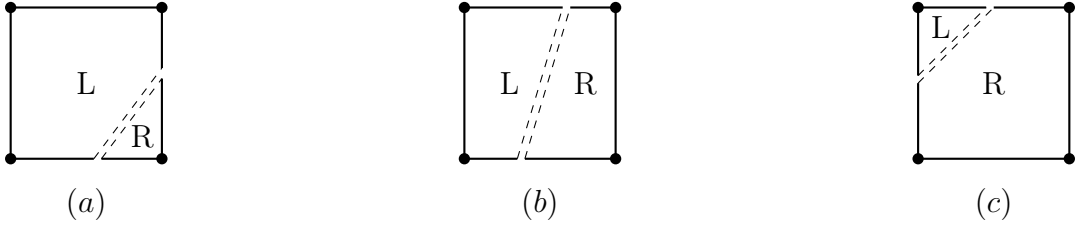
\begin{figure}
\begin{tikzpicture}[scale=2,
    dot/.style={circle, fill=black, inner sep=1.5pt},
    dashedline/.style={dashed}]

\begin{scope}[xshift=0cm]
\draw[thick] (1,0.6) -- (1,1) -- (0,1) -- (0,0) -- (0.55,0);
\draw[thick] (0.6,0) -- (1,0) -- (1,0.53);
\draw[dashedline] (0.55,0) -- (1,0.6);
\draw[dashedline] (0.6,0) -- (1,0.53);
\node[dot] at (0,0) {};
\node[dot] at (1,0) {};
\node[dot] at (1,1) {};
\node[dot] at (0,1) {};
\node at (0.5,0.5) {L};
\node at (0.87,0.15) {R};
\node at (0.5,-0.3) {$(a)$};
\end{scope}

\begin{scope}[xshift=3cm]
\draw[thick] (0.65,1) -- (0,1) -- (0,0) -- (0.35,0);
\draw[thick] (0.4,0) -- (1,0) -- (1,1) -- (0.7,1);
\draw[dashedline] (0.35,0) -- (0.65,1);
\draw[dashedline] (0.4,0) -- (0.7,1);
\node[dot] at (0,0) {};
\node[dot] at (1,0) {};
\node[dot] at (1,1) {};
\node[dot] at (0,1) {}; 
\node at (0.3,0.5) {L};
\node at (0.8,0.5) {R};
\node at (0.5,-0.3) {$(b)$};
\end{scope}

\begin{scope}[xshift=6cm]
\draw[thick] (0,0.5) -- (0,0) -- (1,0) -- (1,1) -- (0.5,1);
\draw[thick] (0,0.55) -- (0,1) -- (0.45,1);
\draw[dashedline] (0,0.5) -- (0.5,1);
\draw[dashedline] (0,0.55) -- (0.45,1);
\node[dot] at (0,0) {};
\node[dot] at (1,0) {};
\node[dot] at (1,1) {};
\node[dot] at (0,1) {};
\node at (0.15,0.87) {L};
\node at (0.5,0.5) {R};
\node at (0.5,-0.3) {$(c)$};
\end{scope}

\end{tikzpicture}
\caption{Some possible orbit spaces of the solution \eqref{ASDCTintro}, cf. also Fig. \ref{FigureSDECT}.}
\label{fig:diagrams}
\end{figure}
\end{center}

\subsection*{Overview}
In \S\ref{sec:THE} we describe the local construction of generic Hermitian-Einstein (non-K\"ahler) toric metrics. In \S\ref{sec:MTE} we solve the Toda equation for the Chen-Teo ansatz, we obtain a constraint, and analyse the solutions; this proves Theorem \ref{thm:solutions}. In \S\ref{sec:ASDE} we study the ASD Einstein case in detail, and we prove Theorem \ref{thm:ASDCT}. In \S\ref{sec:GI} we analyse the construction of gravitational instantons, proving Theorem \ref{thm:instantons}. We conclude in \S\ref{sec:discussion} with a discussion of Delzant polytopes.

\section{Hermitian-Einstein 4-manifolds with two commuting isometries}
\label{sec:THE}

\subsection{Local description}
Here we obtain the local form of a generic toric, Hermitian-Einstein 4-manifold $(M,g,J)$. We assume that $(M,g)$ is not K\"ahler with respect to $J$.

The Hermitian-Einstein condition implies that $(M,g)$ is strictly conformally K\"ahler. The K\"ahler metric is $g_{K}=Z^{-2}g$ where $Z$ is a non-vanishing scalar field. The vector field defined by $\xi_{b}=J^{a}{}_{b}\partial_{a}Z$ is a holomorphic Killing field. We also assume that there is a second holomorphic Killing vector $\eta^{a}$ commuting with $\xi^{a}$ and which is also a Killing vector of $g_{K}$. 

We will first deduce the form of the metric and of the fundamental 2-form $\kappa=g(J\cdot,\cdot)$: there exist local coordinates $(\psi,X,Y,\tilde{X})$ such that they are given by
\begin{align}
g ={}& g_{\psi\psi}\d\psi^2 + 2g_{\psi Y}\d\psi\d Y + g_{YY}\d Y^2
 + g_{\psi\psi}\d\tilde{X}^2 + 2g_{\psi Y}\d\tilde{X}\d{X}+ g_{YY}\d{X}^2, 
 \label{metric0} \\
\kappa ={}& g_{\psi\psi}\d{\psi}\wedge\d\tilde{X} + g_{\psi Y}(\d{\psi}\wedge\d{X}+\d Y \wedge\d\tilde{X}) + g_{YY}\d Y \wedge\d{X}. \label{fund2form0}
\end{align}
To show this we proceed as follows. Let $\tilde{V}^{a}=J^{a}{}_{b}\xi^{b}$ and $V^{a}=J^{a}{}_{b}\eta^{b}$. Integrability of $J$ is equivalent to the vanishing of the Nijenhuis tensor $N_{J}$. In particular, we have $N_{J}(\tilde{V},V)=0$, which implies, after a short calculation, that $\tilde{V}$ and $V$ commute, $[\tilde{V},V]=0$. Thus there are local coordinates $\tilde{X},X$ defined by $\tilde{V}=\frac{\partial}{\partial\tilde{X}}$ and $V=\frac{\partial}{\partial{X}}$. Introducing also Killing coordinates $\psi,Y$ by $\xi^{a}\partial_{a}=\partial_{\psi}$, $\eta^{a}\partial_{a}=\partial_{Y}$, we have a local coordinate system $(\psi,X,Y,\tilde{X})$. 

Since $\eta^a$ is a Killing vector of both $g$ and $g_{K}=Z^{-2}g$, we have $\pounds_{\eta}Z=0$ and $\pounds_{\eta}\kappa_{K}=0$, where $
\kappa_{K}=Z^{-2}\kappa$ is the K\"ahler form. Using that the latter is closed, from Cartan's formula for the Lie derivative we deduce that there is locally a scalar field $H$ such that $\eta_{b}=Z^2J^{a}{}_{b}\partial_{a}H$. Furthermore, $\pounds_{\eta}Z=0$ implies $\kappa^{ab}\partial_{a}Z\partial_{b}H=0$, which gives $\kappa(\xi,\eta)=0$.

Using the above definitions and information, notice that $g_{\psi\psi}=g(\xi,\xi)=g(\tilde{V},\tilde{V})$, $g_{\psi Y}=g(\xi,\eta)=g(\tilde{V},V)$ and $g_{YY}=g(\eta,\eta)=g(V,V)$. Furthermore, $g(\xi,\tilde{V})=g(\xi,J\xi)=\kappa(\xi,\xi)=0$ and $g(\xi,V)=g(\xi,J\eta)=\kappa(\eta,\xi)=0$, and similarly $g(\eta,V)=0=g(\eta,\tilde{V})$. Thus, eq. \eqref{metric0} follows.

Finally, to compute the complex structure, we simply write \eqref{metric0} as 
\begin{align*}
g ={}& g_{\psi\psi}(\d\psi^2+\d\tilde{X}^2)+2g_{\psi Y}(\d\psi \d Y+\d\tilde{X}\d{X})+g_{YY}(\d Y^2+\d{X}^2) \\
={}& g_{\psi\psi}\d{z}^1\d\bar{z}^1+g_{\psi Y}(\d{z}^1\d\bar{z}^2+\d{z}^2\d\bar{z}^1)+g_{YY}\d{z}^2\d\bar{z}^2
\end{align*}
where in the second line we defined $\d{z}^1=\d\psi+\i\d\tilde{X}$ and $\d{z}^2=\d Y+\i\d{X}$. Thus the metric has the Hermitian expression $g = 2g_{\alpha\bar\beta}\d{z}^{\alpha}\d\bar{z}^{\beta}$, so $z^1=\psi+\i\tilde{X}$ and $z^2=Y+\i X$ are holomorphic coordinates. (Here, $g_{\alpha\bar\beta}=g(\frac{\partial}{\partial{z}^{\alpha}}, \frac{\partial}{\partial\bar{z}^{\beta}})$.) The fundamental 2-form is given by $\kappa=\i g_{\alpha\bar\beta}\d{z}^{\alpha}\wedge\d\bar{z}^{\beta}$. Replacing the above expressions, the formula \eqref{fund2form0} follows.

We now turn to the Einstein condition. Starting from \eqref{metric0}-\eqref{fund2form0}, one can show that these formulae are equivalent to the following:
\begin{align}
 g ={}& W^{-1}(\d\psi+\omega\d{Y})^2+W[\d{Z}^2+e^{u}(\d{X}^2+\d{Y}^2)], 
 \label{gLNF} \\
 \kappa ={}& (\d\psi+\omega\d{Y})\wedge\d{Z} + W e^{u} \d{X}\wedge\d{Y},
 \label{f2fLNF}
\end{align}
where we are working away from possible fixed points of $\partial_\psi$, and we defined
\begin{align}\label{Todavariables}
W^{-1}:=g_{\psi\psi}, \qquad 
\omega:=\frac{g_{\psi Y}}{g_{\psi\psi}}, \qquad
e^{u}:=g_{\psi\psi}g_{\psi Y}-(g_{YY})^2.
\end{align}
To prove \eqref{gLNF}-\eqref{f2fLNF}, we use the identity $\d{Z}=g_{\psi\psi}\d\tilde{X}+g_{\psi Y}\d{X}$: this follows from the expression $\xi_{b}=J^{a}{}_{b}\partial_{a}Z$, which is equivalent to $\d{Z}=\partial_{\psi}\lrcorner\,\kappa$. Then we simply write $\d\tilde{X}$ in terms of $\d{Z}$ and $\d{X}$. 

In terms of the formulation \eqref{gLNF}, the Einstein condition $R_{ab}=\lambda g_{ab}$ reduces to the modified $SU(\infty)$ Toda equation with an extra symmetry \cite{Tod:2020ual, a23}: 
\begin{align}\label{MTE}
u_{XX} + (e^{u})_{ZZ}
-\frac{6\lambda Z^2}{\lambda Z^3+3\gamma}(e^{u})_{Z} + \frac{12\lambda Z}{\lambda Z^3+3\gamma}e^{u} = 0.
\end{align}
Here, $\gamma$ is a real constant, such that $\gamma=0$ if and only if the Weyl tensor is anti-self-dual \cite{a23}.
Once a solution to \eqref{MTE} is obtained, the function $W$ is found by the formula
\begin{align}\label{WEinstein}
W=\frac{3 Z(1-\frac{Z}{2}u_{Z})}{\lambda Z^3+3\gamma}.
\end{align}
Finally, knowing $u$ and $W$, the remaining unknown $\omega$ can be found by solving the following system:
\begin{subequations}\label{integrability}
\begin{align}
 \partial_{Z}\omega ={}& \partial_{X}W, \\
 \partial_{X}\omega ={}& -Z^2\partial_{Z}(Z^{-2} W e^u),
\end{align}
\end{subequations}
which follows from the integrability of the complex structure and the condition $\d(Z^{-2}\kappa)=0$.

\subsection{Orthogonal coordinates}
\label{sec:orthogonal}

In \S\ref{sec:MTE} we will need to solve the equations \eqref{MTE} and \eqref{integrability} for the ansatz \eqref{CTansatz}-\eqref{CTansatzfunctions}. Here we give a reformulation of the problem that allows to do this.

We introduce two new coordinates $x,y$ by
\begin{align}\label{orthogonal}
\d{Z} = A\d{x}+B\d{y}, \qquad \d{X} = C\d{x} + D\d{y},
\end{align}
where the functions $A(x,y),B(x,y),C(x,y),D(x,y)$ satisfy $\partial_{y}A=\partial_{x}B$ and $\partial_{y}C=\partial_{x}D$. The dual vector fields are
\begin{equation}\label{dualVF}
\begin{aligned}
\partial_{Z} = \frac{1}{\Delta}\left(D\partial_{x}-C\partial_{y}\right) , \qquad
\partial_{X} = \frac{1}{\Delta}\left(-B\partial_{x}+A\partial_{y}\right),
\end{aligned}
\end{equation}
where $\Delta:=AD-BC\neq0$. Since we are working effectively in two dimensions, it is always possible to require $\partial_{x},\partial_{y}$ to be orthogonal. Computing then $\d Z^2+e^{u}\d X^2$ in terms of $\d x, \d y$, and requiring $g(\partial_{x},\partial_{y})=0$, we obtain
\begin{align}\label{eu}
 e^{u} = \frac{-AB}{CD}.
\end{align}
The metric becomes
\begin{align}\label{metric2}
g = \frac{1}{W}(\d\psi+\omega\d{Y})^2 + We^{u}\d{Y}^2 
+W\Delta\left(\frac{A}{D}\d{x}^2 - \frac{B}{C}\d{y}^2 \right).
\end{align}

We now reformulate the modified Toda equation \eqref{MTE} in terms of $A,B,C,D$: one can show that \eqref{MTE} is equivalent to
\begin{align}\label{MTE2}
\left\{(\lambda Z^3+3\gamma)\left[\partial_{y}\left(\frac{C}{B}\partial_{y}\right) - \partial_{x}\left(\frac{D}{A}\partial_{x}\right)\right] - 6\lambda Z^2 (C\partial_{y}-D\partial_{x})
 - 12\lambda Z \Delta \right\} \frac{AB}{CD} = 0 
\end{align}
where $Z$ is understood as a function $Z(x,y)$. To prove this, we follow \cite{a23}  and introduce two auxiliary variables $a,b$ by
\begin{align}\label{auxvar}
 u_{X} = a_{Z}, \qquad (e^u)_{Z} = -a_{X}+b.
\end{align}
Then equation \eqref{MTE} is satisfied if and only if $b_{Z} = -4\lambda W e^{u}$.
Using \eqref{dualVF}, the two equations in \eqref{auxvar} lead respectively to
\begin{align*}
-B\partial_{x}u + A\partial_{y}u = D\partial_{x}a - C\partial_{y}a, \qquad
\frac{AB}{C}\partial_{x}u - \frac{AB}{D}\partial_{y}u = 
-B\partial_{x}a+A\partial_{y}a - \Delta\, b.
\end{align*}
Writing these equations in matrix form, it is straightforward to find an expression for $(\partial_{x}a,\partial_{y}a)$ in terms of $(\partial_{x}u,\partial_{y}u)$ and $b$:
\begin{align*}
 \partial_{x} a = \frac{A}{D}\partial_{y}u + Cb, \qquad
 \partial_{y} a = \frac{B}{C}\partial_{x}u + Db.
\end{align*}
Now, using $\partial_{y}C = \partial_{x}D$, the expression for $\partial_{Z}$ in \eqref{dualVF}, the condition $b_{Z} = -4\lambda W e^{u}$, and the equality $\partial_{y}\partial_{x} a = \partial_{x}\partial_{y} a$, we obtain
\begin{align}\label{intcondu}
 \partial_{y}\left( \frac{A}{D}\partial_{y}u \right) - \partial_{x}\left( \frac{B}{C}\partial_{x}u \right) = -4\lambda \Delta W e^{u}.
\end{align}
Finally, replacing $e^{u}=-\frac{AB}{CD}$, $\Delta=AD-BC$, and the expression \eqref{MTE} for $W$, into equation \eqref{intcondu}, we arrive at \eqref{MTE2}.

The system \eqref{integrability} can also be conveniently rewritten in terms of $x,y$:
\begin{equation}\label{integrability2}
\begin{aligned}
\partial_{x}\omega ={}& \frac{A}{D}\left[\partial_{y}W - \frac{2B}{Z}\left(1-\frac{Z}{2}u_{Z} \right)W \right], \\
\partial_{y}\omega ={}& \frac{B}{C}\left[\partial_{x}W - \frac{2A}{Z}\left(1-\frac{Z}{2}u_{Z} \right)W \right].
\end{aligned}
\end{equation}
This follows from rewriting \eqref{integrability} as
\begin{align*}
 \partial_{Z}\omega = \partial_{X}W, \qquad
 \partial_{X}\omega = -e^{u}\left[ \partial_{Z}W - \frac{2}{Z}\left(1-\frac{Z}{2}u_Z\right)W\right],
\end{align*}
and, using \eqref{dualVF} and $e^{u}=-\frac{AB}{CD}$, these equations imply
\begin{align*}
D\partial_{x}\omega - C\partial_{y}\omega ={}& -B\partial_{x}W + A\partial_{y}W, \\
-B\partial_{x}\omega + A\partial_{y}\omega ={}& 
\frac{AB}{C}\partial_{x}W - \frac{AB}{D}\partial_{y}W - \Delta\frac{AB}{CD}\frac{2}{Z}\left(1-\frac{Z}{2}u_Z\right)W.
\end{align*}
We can write these equations in matrix form, then it is straightforward to find an expression for $(\partial_{x}\omega,\partial_{y}\omega)$ in terms of $(\partial_{x}W,\partial_{y}W)$ and $W$, which gives \eqref{integrability2}.

\section{The modified Toda equation and its solutions}
\label{sec:MTE}

\subsection{The ansatz}
\label{sec:ansatz}
The formulation of the Einstein condition in section \ref{sec:THE} is generic, applying to any Hermitian-Einstein metric with two commuting Killing fields. In order to apply this approach to the construction of the Chen-Teo metric with $\lambda$, we need to specialise this formulation to obtain an ansatz like \eqref{CTansatz}. 
This can be done in two steps. First, set
\begin{align}\label{AnsatzZ}
Z(x,y) \equiv \frac{\nu x+y }{x-y},
\end{align}
where $\nu$ is an arbitrary parameter. The functions $A$ and $B$ in \eqref{orthogonal} are then
\begin{align}\label{ABansatz}
A = \frac{-(1+\nu)y}{(x-y)^2}, \qquad 
B = \frac{(1+\nu)x}{(x-y)^2}.
\end{align}
Second, recall that the functions $C$ and $D$ in \eqref{orthogonal} must satisfy $\partial_{y}C=\partial_{x}D$. We then choose the simplest solution to this equation:
\begin{align}\label{AnsatzX}
C(x,y) = C(x), \qquad D(x,y) = D(y),
\end{align}
where $C(x)$ and $D(y)$ remain undetermined. Now we simply reparametrise these functions as
\begin{align}\label{CDansatz}
\mathcal{X}(x):=\frac{-(1+\nu)x}{C(x)}, \qquad
\mathcal{Y}(y):=\frac{(1+\nu)y}{D(y)}.
\end{align}
Replacing \eqref{ABansatz} and \eqref{CDansatz} into \eqref{metric2}, we get
\begin{align}\label{metricCT}
g = \frac{(F\d\tau+G\d\phi)^2}{(x-y)FH}+\frac{kH}{(x-y)^3}\left(\frac{\d{x}^2}{\X} - \frac{\d{y}^2}{\Y} - \frac{\X\Y}{kF}\d\phi^2\right),
\end{align}
where 
\begin{align}\label{CTvariables}
F = y^2\mathcal{X} - x^2 \mathcal{Y}, \quad
H = \frac{(\nu+1)^2FW}{k(x-y)}, \quad
G = \frac{(\nu+1)^2F\omega}{k}, \quad
\tau = \frac{(\nu+1)\psi}{\sqrt{k}}, \quad
\phi = \frac{\sqrt{k}\,Y}{(\nu+1)}.
\end{align}
We have thus recovered the Chen-Teo ansatz \eqref{CTansatz} from the assumptions \eqref{AnsatzZ} and \eqref{AnsatzX} \footnote{We note that a natural generalisation of the ansatz \eqref{AnsatzZ} is to consider a generic rational bilinear function of $(x,y)$. However, since we are free to make transformations $x\to x'(x)$, $y\to y'(y)$, this would reduce to \eqref{AnsatzZ} after M\"obius transformations.}. Note that in our current formulation we artificially introduced the constant $k$, so it is not a genuine parameter of the metric.

\subsection{A constraint}
\label{sec:solution}

Replacing the expressions \eqref{ABansatz} and \eqref{CDansatz} for $A,B,C,D$ into the reformulation \eqref{MTE2} of the modified Toda equation, we find 
\begin{equation}\label{MTECT}
\begin{aligned}
 & \left\{  (\lambda Z^3+3\gamma)\left[-\frac{1}{\mathcal{X}}\partial_{y}\left( (x-y)^2\partial_{y} \right) + \frac{1}{\mathcal{Y}}\partial_{x}\left( (x-y)^2\partial_{x} \right) \right] \right. \\
&\left.  +6\lambda Z^2(1+\nu)\left[ \frac{x}{\mathcal{X}}\partial_{y} + \frac{y}{\mathcal{Y}}\partial_{x} \right] + \frac{12\lambda Z(1+\nu)^2 [ y^2\mathcal{X} - x^2\mathcal{Y}]}{(x-y)^2\mathcal{X}\mathcal{Y}} \right\}\frac{\mathcal{X}\mathcal{Y}}{(x-y)^4} 
= 0
\end{aligned}
\end{equation}
where we used that $\Delta = \frac{-F}{(x-y)^2\mathcal{X}\mathcal{Y}}$, together with the explicit expression for $F$ given in \eqref{CTvariables}. After some computations, one can show that the equation is equivalent to
\begin{align}\label{MTECT2}
S(x,y):=f \, (\ddot{\mathcal{Y}}-\ddot{\mathcal{X}})-2 f_{y}\dot{\mathcal{Y}}+2 f_{x}\dot{\mathcal{X}}+2f_{yy}\mathcal{Y}-2f_{xx}\mathcal{X} = 0,
\end{align}
where 
\begin{align}
f(x,y):=\lambda(\nu x+y)^3 + 3\gamma(x-y)^3.
\end{align}

We can now find $\X$ and $\Y$ by taking derivatives of \eqref{MTECT2}. For example, it is straightforward to check that $\partial_{x}^{4}\partial_{y}^{4}S=0$, whereas 
\begin{align}\label{derivatives43}
\partial_{x}^{4}\partial_{y}^{3}S = 6(3\gamma-\lambda)\X^{(6)}, \qquad
\partial_{x}^{3}\partial_{y}^{4}S = 6(3\gamma+\lambda\nu^3)\Y^{(6)},
\end{align}
where $\X^{(n)}$ and $\Y^{(n)}$ denote the $n$-th derivative of $\X$ and $\Y$ with respect to $x$ and $y$, respectively. Since $S(x,y)\equiv0$, from \eqref{derivatives43} we see that, for arbitrary values of the parameters\footnote{After a lengthy computation, the analysis of subcases like $3\gamma=\lambda$, etc. seems to lead to the same conclusions as the generic case, hence we will not include these details here. \label{footnote:subcases}}, we must have $\X^{(6)}=0$ and $\Y^{(6)}=0$. This implies that $\X$ and $\Y$ are fifth-order polynomials
\begin{subequations}\label{X5Y5}
\begin{align}
\X(x) ={}& a_0+a_1 x+a_2 x^2+a_3 x^3+a_4 x^4 + a_5 x^5 \\
\Y(y) ={}& b_0+b_1 y+b_2 y^2+b_3 y^3+b_4 y^4 + b_5 y^5 
\end{align}
\end{subequations}
for some $a_0,...,a_5,b_0,...b_5$. Replacing into \eqref{MTECT2}, we get $S(x,y)=\sum_{i,j=0}^{6}S_{ij}x^{i}y^{j}$, where the coefficients $S_{ij}$, which involve $a_0,...,a_5,b_0,...b_5, \nu,\lambda,\gamma$, must all vanish independently. It turns out 19 out of all the coefficients $S_{ij}$ are non-zero. Their vanishing imposes restrictions on the parameters $a_0,...,a_5,b_0,...b_5, \nu,\lambda,\gamma$. We find:
\begin{align}
\nonumber a_5(3\gamma+\lambda\nu^3) ={}& 0, 
& a_0(3\gamma-\lambda\nu^2) - b_0(3\gamma-\lambda)={}& 0, \\
\nonumber a_5(3\gamma-\lambda\nu^2)={}& 0, 
& a_0(3\gamma+\lambda\nu^3)-b_0(3\gamma+\lambda\nu)={}& 0, \\
\nonumber a_5(3\gamma+\lambda\nu)={}& 0, 
&a_1(3\gamma+\lambda\nu) - b_1(3\gamma-\lambda)={}& 0, \\
\nonumber b_5(3\gamma-\lambda)={}& 0, 
& a_1(3\gamma+\lambda\nu^3)-b_1(3\gamma-\lambda\nu^2)={}& 0, \\
\nonumber b_5(3\gamma-\lambda\nu^2)={}& 0, 
& a_3(3\gamma-\lambda\nu^2)-b_3(3\gamma+\lambda\nu^3)={}& 0, \\
\nonumber b_5(3\gamma+\lambda\nu)={}& 0, 
& a_3(3\gamma-\lambda)-b_3(3\gamma+\lambda\nu)={}& 0, \\
\nonumber a_5(3\gamma-\lambda)+b_5(3\gamma+\lambda\nu^3)={}& 0, 
& a_4(3\gamma+\lambda\nu)-b_4(3\gamma+\lambda\nu^3)={}& 0, \\
\nonumber (a_2-b_2)(3\gamma-\lambda)={}& 0, 
& a_4(3\gamma-\lambda)-b_4(3\gamma-\lambda\nu^2)={}& 0, \\
\nonumber (a_2-b_2)(3\gamma+\lambda\nu^3)={}& 0, \\
\nonumber (a_2-b_2)(3\gamma-\lambda\nu^2)={}& 0, \\
(a_2-b_2)(3\gamma+\lambda\nu)={}& 0. \label{restrictions}
\end{align}
If $\lambda,\gamma,\nu$ are arbitrary (see footnote \ref{footnote:subcases}), then from the first column in \eqref{restrictions} we see that we must have $a_5=b_5=0$ and $a_2=b_2$. Now the equations in the second column still need to be satisfied. For example, the restrictions for $(a_0,b_0)$ can be expressed in matrix form as:
\begin{align}
\left( \begin{matrix} 
3\gamma-\lambda\nu^2 & - (3\gamma-\lambda) \\
3\gamma+\lambda\nu^3 & - (3\gamma+\lambda\nu)
\end{matrix} \right)
\left( \begin{matrix} a_0 \\ b_0 \end{matrix} \right) = 0.
\end{align}
In order to have $a_0\neq0, b_0\neq0$, the determinant of the matrix on the left must vanish. Similar matrix expressions are obtained for $(a_1,b_1),(a_3,b_3),(a_4,b_4)$ from \eqref{restrictions}. It turns out that the determinant is the same for all four pairs, and it has a simple expression:
\begin{align}
\det = 3\lambda\gamma(\nu-1)(\nu+1)^{2}.
\end{align}
To have $a_i\neq0,b_i\neq0$, $i=0,1,3,4$, we must impose $\det=0$, which gives four possibilities:
\begin{align}\label{solutions}
(i) \,\, \lambda=0, \qquad 
(ii) \,\, \gamma=0, \qquad 
(iii) \,\, \nu=1, \qquad 
(iv) \,\, \nu=-1.
\end{align}
The case $\lambda=0$ corresponds to Ricci-flat, and gives the standard Chen-Teo solution. The cases $\gamma=0$ and $\nu=1$ will be analysed in sections \ref{sec:SDE} and \ref{sec:PD}. The case $\nu=-1$ is not allowed at this stage, due to \eqref{AnsatzZ} (but it can be obtained as a limit of the solution with $\gamma=0$). 

Consider now the case $a_i=0,b_i=0$ for $i=0,1,3,4$, then the equations in the second column in \eqref{restrictions} are satisfied, and we have to deal with the first column. We are interested in the case $\lambda\neq0$. If $b_5\neq0$, then we must have $3\gamma=\lambda=\lambda\nu^2=-\lambda\nu$, which implies $\nu=-1$, so this is not allowed. Thus $b_5=0$. Then a similar argument leads to $a_5=0$. Thus the only remaining possibility is $a_2=b_2$. So the solution is in this case $\X = a_2 x^2$, $\Y = a_2 y^2$. But then the function $F$ in \eqref{CTvariables} vanishes, which is not possible. So this case cannot occur.

\subsection{The solution with $\nu=1$}
\label{sec:PD}

The results of section \ref{sec:solution} indicate that, for the ansatz \eqref{AnsatzZ}-\eqref{AnsatzX}, the only possible solutions to the modified Toda equation correspond to \eqref{solutions}. Here we will analyse the case $\nu=1$. The full geometry is found to be

\begin{equation}\label{CTPD}
\begin{aligned}
g ={}& \frac{(F\d\tau+G\d\phi)^2}{(x-y)FH}+\frac{kH}{(x-y)^3}\left(\frac{\d{x}^2}{\mathcal{X}} - \frac{\d{y}^2}{\mathcal{Y}} - \frac{\mathcal{X}\mathcal{Y}}{kF}\d\phi^2\right), \\
\X ={}& a_0 + a_1 x + a_2 x^2 + a_3 x^3 + a_4 x^4, \\ 
\Y ={}& a_0 +  \tfrac{(3\gamma+\lambda)}{(3\gamma-\lambda)}a_1 y 
+ a_2 y^2 + \tfrac{(3\gamma-\lambda)}{(3\gamma+\lambda)}a_3 y^3 + a_4 y^4, \\ 
F ={}& y^2\X-x^2\Y, \\
H ={}& \frac{3(x^2-y^2)}{k}\left( \frac{a_1}{(3\gamma-\lambda)} -\frac{a_3}{(3\gamma+\lambda)} xy \right), \\ 
G ={}& \frac{-3}{k}\left[\left( \frac{a_1y}{(3\gamma-\lambda)} -\frac{a_3y^3}{(3\gamma+\lambda)}  \right)\X+\left( \frac{a_1x}{(3\gamma-\lambda)} -\frac{a_3x^3}{(3\gamma+\lambda)}  \right)\Y \right], 
\end{aligned}
\end{equation}
where $a_0,..,a_4,k,\lambda$ are real constants, such that $\lambda\neq\pm3\gamma$. 

To prove this, we replace \eqref{X5Y5} into \eqref{MTECT2} and set $\nu=1$. We find that the solution to $S\equiv0$ is 
\begin{align*}
 b_5 = a_5 = 0, \qquad b_0=a_0, 
 \qquad b_1= \tfrac{(3\gamma+\lambda)}{(3\gamma-\lambda)}a_1,
 \qquad b_2=a_2,
 \qquad b_3= \tfrac{(3\gamma-\lambda)}{(3\gamma+\lambda)}a_3,
 \qquad b_4=a_4,
\end{align*}
where we are assuming the generic case $\lambda\neq\pm3\gamma$. This gives the polynomials $\X$ and $\Y$ in \eqref{CTPD}. Knowing $\X$ and $\Y$, we determine $u$ via \eqref{eu}, and we use this to find $W$ via \eqref{WEinstein}. Then we can find $\omega$ by integrating the system \eqref{integrability2}. Replacing into \eqref{CTvariables}, we find the functions $H$ and $G$ as given in \eqref{CTPD}.

We now demonstrate that the metric is isometric to Pleba\'nski-Demia\'nski with $\lambda$. Note first that we have a scaling freedom for the parameters $a_0,...,a_4$, which we can use to set two of the parameters to any non-zero values. We can then choose 
$ a_1 = \frac{\lambda-3\gamma}{12}$, $a_3 = \frac{\lambda+3\gamma}{12}$,
which gives 
\[
H =-\frac{1}{4k}(x^2-y^2)(1+xy), \qquad
G = \frac{1}{4k}\left[(y+y^3)\X+(x+x^3)\Y \right].
\]
We now perform a M\"obius transformation $x = \frac{p+1}{p-1}$, $ y = \frac{q+1}{q-1}$ (which is inspired by \cite[Eq. (4.2)]{CT2}),
and define new parameters
\[
c_{0,4}=a_0+a_2+a_4\mp\frac{\lambda}{6}, \qquad 
c_{1,3} = -4a_0+4a_4\mp\gamma, \qquad
c_2 = 6(a_0+a_4)-2a_2,
\]
together with the polynomials
\[
P = c_0+c_1 p+c_2 p^2+c_3 p^3+c_4 p^4, \qquad 
Q = c_0+\tfrac{\lambda}{3}+c_1 q+c_2 q^2+c_3 q^3+(c_4-\tfrac{\lambda}{3}) q^4.
\]
After a short calculation, the non-Killing part of the metric \eqref{metricCT} becomes
\begin{align}
\frac{kH}{(x-y)^3}\left(\frac{\d x^2}{\X}-\frac{\d y^2}{\Y}\right) = 
\frac{(1-p^2q^2)}{(p-q)^2}\left( \frac{\d p^2}{P}-\frac{\d q^2}{Q} \right).
\label{nonKillingPartPD}
\end{align}
Transforming now from $(\tau,\phi)$ to new Killing coordinates $(\tilde\tau,\tilde\phi)$ given by $\tau = \tfrac{1}{\sqrt{k}} \, (\tilde\tau+\tilde\phi)$, $\phi = 2\sqrt{k} \, (\tilde\tau-\tilde\phi)$,
the Killing part of the metric \eqref{metricCT} becomes
\begin{align}
\frac{(F\d\tau+G\d\phi)^2}{(x-y)FH}-\frac{H\X\Y \, \d\phi^2}{(x-y)^3F} 
= \frac{(q^4P-Q)\d\tilde\tau^2+(p^2Q-q^2P)\d\tilde\tau\d\tilde\phi+(P-p^4Q)\d\tilde\phi^2}{(p-q)^2(1-p^2q^2)}.
\label{KillingPartPD}
\end{align}
Using then \eqref{metricCT} together with \eqref{nonKillingPartPD} and \eqref{KillingPartPD}, we see that the metric \eqref{metricCT} is the Pleba\'nski-Demia\'nski solution with cosmological constant $\lambda$ (see e.g. \cite[Prop. 4.3]{a23}).

\subsection{The solution with $\gamma=0$}
\label{sec:SDE}

From section \ref{sec:solution}, $\gamma=0$ is the other non-trivial solution to the restrictions \eqref{restrictions}. As mentioned in section \ref{sec:THE}, the case $\gamma=0$ in the modified Toda equation \eqref{MTE} corresponds to an anti-self-dual Weyl tensor. The local solution is found to be 
\begin{equation}\label{SDECT}
\begin{aligned}
g ={}& \frac{(F\d\tau+G\d\phi)^2}{(x-y)FH}+\frac{kH}{(x-y)^3}\left(\frac{\d{x}^2}{\mathcal{X}} - \frac{\d{y}^2}{\mathcal{Y}} - \frac{\mathcal{X}\mathcal{Y}}{kF}\d\phi^2\right), \\
\X ={}& a_0 + a_1 x + a_2 x^2 + a_3 x^3 + a_4 x^4, \\
\Y ={}& a_0\nu^2 - a_1\nu y + a_2 y^2 - \frac{a_3}{\nu} y^3 + \frac{a_4}{\nu^2} y^4, \\
 F ={}& y^2 \X - x^2 \Y,  \\
H ={}& \frac{3(\nu+1)}{2k\lambda}(x-y)\left[ 2(\nu-1)(a_0-\frac{a_4}{\nu^2}x^2y^2) - (x+y)(a_1+\frac{a_3}{\nu}xy)\right], \\
G ={}& \frac{3(\nu+1)}{2k \lambda \nu}\left[ 
\left(a_3 y^3-2(\nu-1)a_2 y^2 - 2a_0(\nu-1)\nu^2+a_1\nu(2\nu-1)y \right) \X \right. \\
& \qquad \left. + \left(2a_2(\nu-1)x^2+2a_4(\nu-1)x^4+a_1\nu x+a_3(2\nu-1)x^3 \right)\Y \right]
\end{aligned}
\end{equation}
where $a_0,...,a_4,\nu,k,\lambda$ are real parameters. 

To show this, we replace \eqref{X5Y5} into \eqref{MTECT2} and set $\gamma=0$; then we find that the solution to $S\equiv0$ is 
\begin{align}\label{zerogamma}
a_5=b_5=0, \qquad b_0=\nu^2 a_0, \qquad b_1=-\nu a_1, \qquad b_2=a_2, \qquad b_3=-\frac{a_3}{\nu}, \qquad b_4=\frac{a_4}{\nu^2}.
\end{align}
Note that this condition implies only that the Weyl tensor is ASD; this solution was actually obtained in \cite[Theorem 5.3]{a23}. The Einstein condition has not yet been imposed: this is obtained after setting the function $W$ to be equal to \eqref{WEinstein}, where $u$ is given by \eqref{eu}. 
Integrating \eqref{integrability2}, we find $\omega$. Once $W$ and $\omega$ are known, the functions $H$ and $G$ in \eqref{metricCT} are obtained from \eqref{CTvariables}. This procedure then leads to the metric \eqref{SDECT}.

\section{The ASD Einstein case}
\label{sec:ASDE}

In this section we will study the solution \eqref{SDECT} in detail, proving the properties summarised in Theorem \ref{thm:ASDCT}. We will be interested in the case $\lambda<0$. 
We postpone the discussion of smoothness (i.e., of gravitational instantons) to section \ref{sec:GI}.

\subsection{Conformal boundary}
\label{sec:boundary}

We first want to show that the solution is asymptotically (locally) hyperbolic. Alternative names for this condition are Poincar\'e-Einstein or conformally compact. We recall that a Riemannian manifold $(M,g)$ is said to satisfy this condition if there exists a compact manifold $\overline{M}$ with boundary $S=\partial\overline{M}$ such that $M=\overline{M}\setminus S$, and a boundary defining function $\Omega:\overline{M}\to\mathbb{R}^{+}$ (i.e., $\Omega|_{S}=0$, $\d\Omega|_{S}\neq0$, $\Omega|_{M}\neq0$) such that $\Omega^2g$ is a non-degenerate metric on $M\cup S=\overline{M}$. The set $\{\Omega=0\}$ is referred to as the conformal boundary (or conformal infinity) of $g$. Near the boundary, the Fefferman-Graham expansion gives 
\begin{align}\label{FeffermanGraham}
g = \frac{\d\Omega^2+h}{\Omega^2} + O(\Omega^{-1})
\end{align}
for some tensor field $h$, which is independent of $\Omega$. The Riemannian 3-manifold $(S,h)$ is called the boundary geometry. The metric $\Omega^2g$ is said to give a conformal compactification of $g$.

To show that the solution \eqref{SDECT} is asymptotically locally hyperbolic, we can use the Toda framework from section \ref{sec:THE}; cf. also \cite{ALR26}. Recall that in this formulation, the metric is given by \eqref{gLNF}, with $u,W$ satisfying \eqref{MTE}-\eqref{WEinstein} for $\gamma=0$. It is convenient to define $\hat{Z}:=\frac{-1}{Z}$, $e^{\hat{u}}:=\frac{e^{u}}{Z^4}$, $\hat{W}:=Z^2W$. Then a short calculation shows that the metric can be written as 
\begin{align}\label{Todahat}
g = \frac{1}{\hat{Z}^2}\left[\hat{W}^{-1}(\d\psi+\omega\d Y)^2+\hat{W}(\d\hat{Z}^2+e^{\hat{u}}(\d X^2+\d Y^2)) \right],
\end{align}
and $\hat{W}=-\frac{3}{\lambda}(1-\frac{\hat{Z}}{2}\hat{u}_{\hat{Z}})$. We will be interested in the case $\lambda<0$, so for the discussion that follows, it is convenient to write $\ell^2:=-3/\lambda$.

For our solution, we have $\hat{Z} = -\frac{(x-y)}{\nu x+y}$, so the zero set corresponds to $x-y=0$. Using the explicit form of the functions $\hat{u},\hat{W},\omega$, we can then write, near $\hat{Z}=0$, the expansions $\hat{u}=\hat{u}_{0}+O(\hat{Z})$, $\hat{W}=\ell^2(1+O(\hat{Z}))$, $\omega=\omega_{0}+O(\hat{Z})$, where $\hat{u}_{0},\omega_{0}$ are independent of $\hat{Z}$. Replacing into \eqref{Todahat}, we get
\begin{align}
g =\frac{1}{\hat{Z}^2}\left[ \ell^2\d\hat{Z}^2+\frac{1}{\ell^2}(\d\psi+\omega_{0}\d Y)^2+\ell^2 e^{\hat{u}_{0}}(\d X^2+\d Y^2) \right] + O(\hat{Z}^{-1}).
\end{align}
Comparing now to \eqref{FeffermanGraham}, we then see that we can interpret $\hat{Z}$ as a boundary defining function, and the set $S=\{\hat{Z}=0\}$ as a conformal boundary. The tensor field $\hat{Z}^2g$ extends through the boundary to a non-degenerate metric (which will be studied further in section \ref{sec:kahler} below). The solution \eqref{SDECT} is then asymptotically hyperbolic, with a conformal boundary at $x-y=0$. To determine the topology of the boundary, note first that the vector fields $\partial_{\psi},\partial_{Y}$ generate an isometric action of the 2-torus $T^2$ on $M$, which extends to $S$. In \S\ref{sec:rods} we will see that this action on $S$ has exactly two singular orbits. It then follows that $S$ must be a lens space. This shows point \eqref{item:ALH} in Theorem \ref{thm:ASDCT}.

\subsection{Global structure}
\label{sec:global}

Having determined that the solution is asymptotically hyperbolic, we can now analyse the global structure in more detail. In order to do this, we need to specify the domain of the coordinates $(x,y)$. We first note that the polynomials $\X$ and $\Y$ in \eqref{SDECT} are slightly different. It is convenient to redefine the variables so that the polynomials are the same. This can be done as follows:
\begin{align}\label{redefcoord}
y \to -\nu y, \qquad \tau \to -\frac{\tau}{\nu}, \qquad G \to \nu G, \qquad
\Y \to \nu^2 \Y, \qquad F \to \nu^2 F.
\end{align}
After these redefinitions, the geometry \eqref{SDECT} becomes
\begin{equation}\label{metricCT2}
\begin{aligned}
g ={}& \frac{({F}\d\tau+{G}\d\phi)^2}{(x+\nu{y}){F}H}
+\frac{kH}{(x+\nu{y})^3}\left(\frac{\d{x}^2}{\X} - \frac{\d{y}^2}{{\Y}} - \frac{\X{\Y}}{k{F}}\d\phi^2\right), \\
\X ={}& a_0 + a_1 x + a_2 x^2 + a_3 x^3 + a_4 x^4, \\
{\Y} ={}& a_0 + a_1{y} + a_2 y^2 +a_3{y}^3 + a_4 {y}^4, \\
 {F} ={}& {y}^2 \X - x^2 {\Y},  \\
H ={}& \frac{3(\nu+1)}{2k\lambda}(x+\nu{y})\left[ 2(\nu-1)(a_0 - a_4 x^2{y}^2) - (x-\nu{y})(a_1 - a_3 x {y})\right], \\
{G} ={}& \frac{3(\nu+1)}{2k \lambda \nu}\left[ 
\left(-a_3 \nu{y}^3-2(\nu-1)a_2 {y}^2 - 2a_0(\nu-1)-a_1(2\nu-1){y} \right) \X \right. \\
& \qquad\qquad \left. + \left(2a_2(\nu-1)x^2+2a_4(\nu-1)x^4+a_1\nu x+a_3(2\nu-1)x^3 \right){\Y} \right],
\end{aligned}
\end{equation}
which is the form of the metric we gave in Theorem \ref{thm:ASDCT}.

We can now describe the global structure in terms of $x$ and the new coordinate $y$. This is similar to description of the Ricci-flat Chen-Teo metric \cite{CT2}, but with the difference that there is now a conformal boundary that can intersect the coordinate domain. Note that after the redefinition \eqref{redefcoord} of $y$, the boundary now corresponds to the line $x+\nu y=0$.

We have $\X(x)=P(x)$ and $\Y(y)=P(y)$, where $P(z)=a_0 + a_1 z + a_2 z^2 + a_3 z^3 + a_4 z^4$. We will assume that the polynomial $P$ has four real roots $r_1,...,r_4$, so that 
\begin{align}\label{roots}
P(z)=a_4(z-r_1)(z-r_2)(z-r_3)(z-r_4). 
\end{align}
Then $x$ must lie between two adjacent roots, and similarly for $y$. Moreover, to ensure Riemannian signature, $x$ and $y$ must lie in different ranges. For concreteness, we will take the box $B:=\{r_2<x<r_3\}\cap\{r_1<y<r_2\}$ in the $(x,y)$-plane. We illustrate the situation in Fig. \ref{FigureSDECT}. 

As mentioned, the geometry has a conformal boundary along the line $x+\nu y=0$. Depending on the value of $\nu$, this line may or may not pass through the box. If the line does not intersect the box, the manifold is compact. The case $\lambda<0$ is then ruled out since the metric has Killing fields. If $\lambda>0$, then we have an ASD Einstein metric on a compact manifold with positive Einstein constant. The only smooth such possibilities are $S^4$ and $\mathbb{CP}^2$. Therefore, to have a more interesting smooth solution we need the line to intersect the box, which requires $\lambda<0$, and is what we will assume from now on.

The domain is then the intersection of the box $B=\{r_2<x<r_3\}\cap\{r_1<y<r_2\}$ with one of the two regions $x+\nu y\lessgtr0$. However, the metric is well-defined on both regions $B\cap\{x+\nu y>0\}$ and $B\cap\{x+\nu y<0\}$; thus, we have {\it two} solutions, each defined on either side of conformal infinity. The domain of each solution is a polygon with either 3, 4, or 5 edges, where one of the edges is the boundary, see Fig. \ref{FigureSDECT}. The number of edges depends on the slope of the line $x+\nu y=0$, i.e. on the value of $\nu$. This implies that the two solutions will in general have different topologies. This shows point \eqref{item:ASDE2} in Theorem \ref{thm:ASDCT}.

\begin{center}
\begin{figure}
\begin{tikzpicture}[scale=1.1]

\draw[->, thick] (-4,0) -- (1,0) node[right] {$x$};
\draw[->, thick] (0,-4) -- (0,1) node[above] {$y$};

\def\rone{-2.5}
\def\rtwo{-1.2}
\def\rthree{1.2}

\def\xleft{-1.5}
\def\xright{-0.3}
\def\ybottom{-3.4}
\def\ytop{-1.5}

\draw[dashed] (\xleft,0) -- (\xleft,\ybottom);
\draw[dashed] (\xright,0) -- (\xright,\ybottom);
\draw[dashed] (0,\ybottom) -- (\xright,\ybottom);
\draw[dashed] (0,\ytop) -- (\xleft,\ytop);

\draw[thick]
  (\xleft,\ybottom) rectangle (\xright,\ytop);

\draw[blue, thick, domain=-2.2:0.5] plot (\x,{\x/0.68});
\node[blue] at (-3,-3) {${y}=-\frac{x}{\nu_3}$};

\draw[orange, thick, domain=-1.4:0.4] plot (\x,{\x/0.37});
\node[orange] at (-1.8,-4) {${y}=-\frac{x}{\nu_2}$};

\draw[red, thick, domain=-0.6:0.15] plot (\x,{\x/0.15});
\node[red] at (-0.25,-4.3) {${y}=-\frac{x}{\nu_1}$};

\node[right] at (0,\ybottom) {$r_1$};
\node[right] at (0,\ytop) {$r_2$};
\node[above] at (\xleft,0) {$r_2$};
\node[above] at (\xright,0) {$r_3$};

\end{tikzpicture}
\caption{Possible domains for $(x,y)$, assuming $r_1<r_2<r_3<0$ and $-1<\nu<0$. The different slanted lines $x+\nu_{i} y=0$ correspond to different locations of conformal infinity, depending on the value of $\nu$.}
\label{FigureSDECT}
\end{figure}
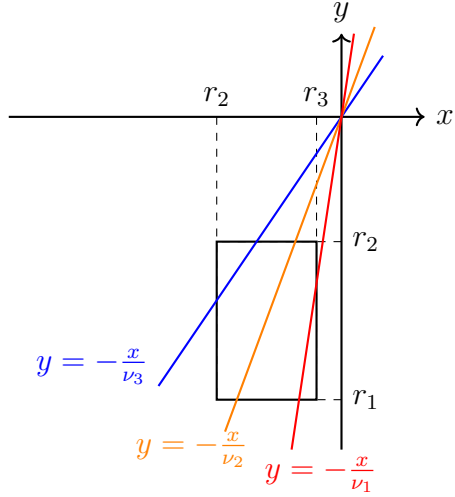
\end{center}

\subsection{K\"ahler geometry and ALE ends}
\label{sec:kahler}

As mentioned in section \ref{sec:THE}, since the metric $g$ is Hermitian-Einstein with $\lambda\neq0$, the metric $g_{K}:=\hat{Z}^2g$ is K\"ahler, with $\hat{Z}=-Z^{-1}$. And since $g$ is ASD, the K\"ahler metric $g_{K}$ is scalar-flat, i.e. it has zero scalar curvature. Being scalar-flat K\"ahler, $g_{K}$ is automatically Einstein-Maxwell \cite{Flaherty}.

We note that, since any Killing vector on an ASD Einstein space defines a conformal K\"ahler structure \cite{Tod2}, and \eqref{metricCT2} has two Killing fields, any linear combination of them defines an independent K\"ahler structure. We will show that $g_{K}$ has the special property of being ALE.

In the new coordinates \eqref{redefcoord}, the conformal factor is 
\begin{align}\label{conformalfactor}
\hat{Z} = -\frac{(x+\nu y)}{\nu(x-y)}.
\end{align}
Note that $|\hat{Z}|\to\infty$ when $x-y\to 0$. For the coordinate domain specified in section \ref{sec:global}, the only point which is also on the line $x-y=0$ is $(x,y)=(r_2,r_2)$. 
Since the conformal infinity line $x+\nu y=0$ cuts through the box, and the two Einstein metrics are defined on opposite sides of this boundary, the solution that does not contain the point $(r_2,r_2)$ is globally conformally K\"ahler, while the other solution is not. 

From the discussion in \S\ref{sec:boundary}, note that the K\"ahler metric $g_K=\hat{Z}^{2}g$ gives a conformal compactification of $g$. It extends beyond the conformal boundary, so it is defined on the whole box $B$. We will now show that the point $(r_2,r_2)$ is an ALE end of $g_{K}$.
To see this, we introduce a new coordinate system $(\Psi,\Phi,R,\Theta)$ as follows:
\begin{equation}
\begin{aligned}
 \tau ={}& \frac{3(\nu^2-1)}{2\lambda\sqrt{k}}(\Psi+b_1\Phi), \\
 \phi ={}& \frac{3(\nu^2-1)}{2\lambda\sqrt{k}}b_2\Phi, \\
 x ={}& r_2 + \frac{c}{2R^2}\cos^2(\Theta/2), \\ 
 y ={}& r_2 - \frac{c}{2R^2}\sin^2(\Theta/2),
\end{aligned}
\end{equation}
where 
\begin{equation}
\begin{aligned}
b_1 ={}& \frac{(3-\nu)r_2(r_1r_2+r_1r_3+r_1r_4+r_2r_3+r_2r_4+r_3r_4)+(1+\nu)(r_2^3+r_1r_3r_4)}{(\nu-1)(r_1 - r_2) (r_2 - r_3) (r_2 - r_4) }, \\
b_2 ={}& -\frac{4k r_2\lambda}{3a_4(\nu^2-1)(r_1 - r_2) (r_2 - r_3) (r_2 - r_4)}, \\
c ={}& \frac{12 r_2(1-\nu^2)}{\lambda\nu^2}.
\end{aligned}
\end{equation}
Notice that $(x,y)\to(r_2,r_2)$ corresponds to $R\to\infty$. Then, one can check that as $R\to\infty$,
\begin{align}\label{ALEKahler}
g_{K} = \mathring{g} + O(R^{-2}),
\end{align}
where 
\begin{align}
\mathring{g}=\d{R}^2+\frac{R^2}{4}\left( (\d\Psi+\cos\Theta\d\Phi)^2+\d\Theta^2+\sin^2\Theta\d\Phi^2 \right),
\end{align}
and the falloff in \eqref{ALEKahler} is estimated with respect to the asymptotic metric $\mathring{g}$. Therefore, with the appropriate identifications of the angles $\Psi,\Phi$, we see that the K\"ahler metric $g_{K}$ is ALE of order $2$ (see e.g. \cite{BKN} for the definition of these terms). This proves point \eqref{item:SFK} in Theorem \ref{thm:ASDCT}.

\subsection{Multipoles}
\label{sec:multipoles}

As mentioned in the introduction, any ASD toric Einstein metric can be described, locally, via the Calderbank-Pedersen construction \cite{CP}. The details can be found in \cite{CP} and we will not review them here, but the key point is that any such metric can be generated in terms of a function $\mathcal{F}$ of two variables $\rho,\zeta$, satisfying 
\begin{align}\label{eigenfunction}
\mathcal{F}_{\rho\rho}+\mathcal{F}_{\zeta\zeta} = \frac{3}{4\rho^2}\mathcal{F}.
\end{align}
That is, $\mathcal{F}$ is an eigenfunction of the Laplacian on the hyperbolic 2-plane $\rho^{-2}(\d\rho^2+\d\zeta^2)$, with eigenvalue $3/4$. 

A class of solutions to \eqref{eigenfunction} that arises naturally in \cite{CP} are the so-called $m$-pole solutions, which are those of the form \eqref{mpole}. Here, we will show that the metric \eqref{metricCT2} is a 4-pole solution, that is, it can be generated by a function of the form \eqref{mpole} with $m=4$. 

It will also be instructive, especially for section \ref{sec:GI}, to describe the solution \eqref{metricCT2} in terms of Tod's formulation of ASD toric Einstein metrics \cite[Section 4.1]{Tod:2020ual}. In this formulation, any such metric can be generated by a solution $V$ to the axisymmetric Laplace equation in $\mathbb{R}^3$, 
\begin{align}\label{Laplace}
V_{\zeta\zeta} + \rho^{-1}(\rho V_{\rho})_{\rho}=0.
\end{align}
Given a solution to \eqref{Laplace}, the metric is obtained as follows (cf. \cite{Tod:2020ual} for details):
\begin{align}\label{Tod1}
g = {}&\frac{1}{\hat{Z}^2}\left[\hat{W}^{-1}(\d\psi+\omega\d Y)^2+\hat{W}\rho^2\d{Y}^2+e^{2\hat\mu}(\d\rho^2+\d\zeta^2) \right], 
\end{align}
where
\begin{align}\label{Tod2}
\hat{Z} ={}& \frac{\rho V_{\rho}}{2}, \qquad
\hat{W} = -\frac{3}{\lambda}\left(1+\frac{V_{\rho}V_{\zeta\zeta}}{\rho\tilde\Delta} \right), \qquad
 \omega = -\frac{3}{\lambda}\left(\zeta - \frac{V_{\rho}V_{\rho\zeta}}{\tilde\Delta} \right), \qquad
 e^{2\hat\mu} = \frac{W\rho^2\tilde\Delta}{4},
\end{align}
and $\tilde\Delta=V_{\zeta\zeta}^2+V_{\rho\zeta}^2$. Here, the coordinates $\psi,Y,\hat{Z}$ and the functions $\omega, \hat{W}$ are the same as in \eqref{Todahat}, while the coordinates $\rho,\zeta$ are defined by 
\begin{align}\label{WPcoord}
\rho^2 = e^{\hat{u}}, \qquad \zeta_{X} = \rho\rho_{\hat{Z}}, \qquad \zeta_{\hat{Z}} = -\rho^{-1}\rho_{X}.
\end{align}

To connect the description based on a solution to \eqref{Laplace} to the description based on a solution to \eqref{eigenfunction}, simply define $\mathcal{F}:=\rho^{-1/2}\hat{Z}$. In particular, this means that the $m$-pole solutions are those in which the conformal factor that makes the metric K\"ahler is of the form
\begin{align}
\hat{Z} = \sum_{k=0}^{m-1}\sqrt{\alpha_{k}\rho^2+(\alpha_{k}\zeta-\beta_k)^2}
\end{align}
for some constants $\alpha_k,\beta_{k}$ (where, compared to \eqref{mpole}, we have for notational convenience redefined the sum from $0$ to $m-1$ instead of $1$ to $m$).

In our case, the conformal factor is given by \eqref{conformalfactor}, and we need to find the coordinates $\rho,\zeta$. This is done via \eqref{WPcoord}, which, after a short computation, gives
\begin{align}\label{WPcoord2}
\rho = \frac{\sqrt{-\X\Y}}{\nu^2(x-y)^2}, \qquad 
\zeta=-\frac{[2 a_0 + a_1 (x + y) + x y (2 a_2 + 2 a_4 x y + a_3 (x + y))]}{
 2 \nu^2(x - y)^2}.
\end{align}
Using then the explicit expressions \eqref{conformalfactor} and \eqref{WPcoord2}, one can check that the following identity is true:
\begin{align}\label{4pole}
\hat{Z} = A + \alpha_1\sqrt{\rho^2+(\zeta-\zeta_1)^2}+\alpha_2\sqrt{\rho^2+(\zeta-\zeta_2)^2}+\alpha_3\sqrt{\rho^2+(\zeta-\zeta_3)^2}, 
\end{align}
where, in terms of the roots introduced in \eqref{roots}, we have
\begin{align}\label{coeff4pole}
A = \frac{\nu-1}{2\nu}, \quad 
\zeta_1 = \frac{a_ 4(r_ 1r_ 2+r_ 3r_ 4)}{2\nu^2}, \quad 
\zeta_2 = \frac{a_ 4(r_ 1r_ 3+r_ 2r_ 4)}{2\nu^2}, \quad 
\zeta_3 = \frac{a_ 4(r_ 2r_ 3+r_ 1r_ 4)}{2\nu^2}
\end{align}
and 
\begin{equation}\label{coeff4pole2}
\begin{aligned}
\alpha_1 = \frac{(r_ 1r_ 2-r_ 3r_ 4)\nu(1+\nu)}{a_ 4 r_{13}r_{23}r_{14}r_{24}}, \qquad
\alpha_2 = \frac{(r_ 1r_ 3-r_ 2r_ 4)\nu(1+\nu)}{a_ 4 r_{12}r_{23}r_{14}r_{34}}, \qquad
\alpha_3 = \frac{(r_ 2r_ 3-r_ 1r_ 4)\nu(1+\nu)}{a_ 4 r_{12}r_{13}r_{24}r_{34}},
\end{aligned}
\end{equation}
with $r_{ij}:=r_{i}-r_{j}$. This can be checked by expressing the coefficients $a_0,a_1,a_2,a_3$ in terms of the roots:
\begin{equation}\label{coeffroots}
\begin{aligned}
a_0 ={}& a_4 r_1 r_2 r_3 r_4, \\
a_1 ={}& -a_4 (r_1 r_2 r_3 + r_1 r_2 r_4 + r_1 r_3 r_4 + r_2 r_3 r_4), \\
a_2 ={}& a_4 (r_1 r_2 + r_1 r_3 + r_2 r_3 + r_1 r_4 + r_2 r_4 + r_3 r_4), \\
a_3 ={}& -a_4 (r_1 + r_2 + r_3 + r_4).
\end{aligned}
\end{equation}

Therefore, we see that the metric \eqref{metricCT2} is a 4-pole solution in the Calderbank-Pedersen sense \cite{CP}. This proves point \eqref{item:CP} in Theorem \ref{thm:ASDCT}, for $\nu\neq\pm1$. (The limits $\nu\to\pm1$ will be briefly addressed in \S\ref{sec:specialnu} below.)
Furthermore, given the explicit form \eqref{4pole}, the corresponding axisymmetric harmonic function $V$ needed for Tod's description \eqref{Tod1}-\eqref{Tod2} is easy to find, by integrating $\hat{Z}=\frac{1}{2}\rho V_{\rho}$:
\begin{align}\label{VASDCT}
V = A\log\rho^2+\sum_{i=1}^{3}\alpha_i\left[2\sqrt{\rho^2+(\zeta-\zeta_i)^2}-(\zeta-\zeta_i)\log\left(\frac{\sqrt{\rho^2+(\zeta-\zeta_i)^2}+(\zeta-\zeta_i)}{\sqrt{\rho^2+(\zeta-\zeta_i)^2}-(\zeta-\zeta_i)} \right) \right],
\end{align}
up to an (irrelevant) affine function of $\zeta$. The description \eqref{Tod1}-\eqref{Tod2}, with $V$ given by \eqref{VASDCT}, will be instrumental for the construction of gravitational instantons in section \ref{sec:GI}.

When the solution \eqref{metricCT2} is generated via \eqref{VASDCT}, \eqref{Tod1}, \eqref{Tod2}, {\it a priori} it would seem to depend on 7 parameters: $A, \alpha_1, \alpha_2, \alpha_3,\zeta_1,\zeta_2,\zeta_3$. However, as explained in \cite{CP} or in more detail in \cite[\S C.2]{Farquet:2014kma}, the $m$-pole solution has several internal symmetries that reduce the number of independent parameters to $2m-4$, which then gives 4 parameters for \eqref{metricCT2}. This matches the counting in terms of the initial parametrisation, where the parameters were $a_0,...,a_4,\nu$ but two of the $a_i$ coefficients are redundant.

\subsection{The subcases $\nu=\pm1$}
\label{sec:specialnu}

With the form \eqref{metricCT2} of the metric, we can analyse the limit $\nu\to-1$ as follows. Note that \eqref{metricCT2} can be written as $g = -\frac{3(\nu+1)}{2\lambda k}\hat{g}$, where 
\begin{align}\label{gnu-1}
\hat{g} = \frac{({F}\d\hat\tau+\hat{G}\d\phi)^2}{(x+\nu{y}){F}\hat{H}}
+\frac{k\hat{H}}{(x+\nu{y})^3}\left(\frac{\d{x}^2}{\X} - \frac{\d{y}^2}{{\Y}} - \frac{\X{\Y}}{k{F}}\d\phi^2\right)
\end{align}
and 
\begin{align}
\hat\tau=\frac{2\lambda k \tau}{3(\nu+1)}, \qquad \hat{G}=\frac{2\lambda k G}{3(\nu+1)}, \qquad \hat{H}=\frac{-2\lambda k H}{3(\nu+1)}.
\end{align}
A short computation shows that $\hat{H}|_{\nu=-1}=H_{\rm CT}|_{\nu=-1}$ and $\hat{G}|_{\nu=-1}=G_{\rm CT}|_{\nu=-1}-a_2 F$, where $H_{\rm CT}, G_{\rm CT}$ are the corresponding functions of the Ricci-flat Chen-Teo metric \cite[Eq. (2.1)]{CT2}. Therefore (redefining $\hat\tau\to \hat\tau+a_2\phi$), the limit $\nu\to-1$ of \eqref{gnu-1} is the same as the limit $\nu\to-1$ of the Ricci-flat Chen-Teo metric, which is a three-centred ALE Gibbons-Hawking space \cite[\S 4.2]{CT2}

For the case $\nu=1$, it is convenient to use the expression \eqref{SDECT} for the metric. Setting $\nu=1$ in \eqref{SDECT}, it is straightforward to show that the metric reduces to \eqref{CTPD} with $\gamma=0$. Since we showed that \eqref{CTPD} is isometric to the Pleba\'nski-Demia\'nski solution with $\lambda$, and the case $\gamma=0$ is ASD, we conclude that the metric \eqref{metricCT2} with $\nu=1$ reduces to ASD Pleba\'nski-Demia\'nski with $\lambda$. This finalises the proof of point \eqref{item:CP} in Theorem \ref{thm:ASDCT}. Note that, since we see from \eqref{coeff4pole} that $A=0$ for $\nu=1$, it follows that ASD Pleba\'nski-Demia\'nski with $\lambda$ is a 3-pole solution.

\section{Gravitational instantons}
\label{sec:GI}

In this section we prove Theorem \ref{thm:instantons}, i.e., we show the existence of gravitational instantons within the family of ASD Einstein metrics \eqref{ASDCTintro}. Since the family is toric, the problem can be approached by first understanding the rod structure of the solution \cite{CT0}, and then analysing the associated regularity conditions. This framework is closely related to the analysis of Delzant polytopes in toric K\"ahler geometry. Our analysis is based on \cite{CT0, BG, ALR26}.

\subsection{Preliminaries}
\label{sec:rodstructure}

Let $(M,g)$ be a Riemannian 4-manifold $(M,g)$ with an isometric action of the 2-torus $T^2=S^1\times S^1$, cf. \cite{OR}. The rod-structure formalism \cite{CT0} associates with $(M,g)$ certain invariants determined by the degeneration of the $T^2$-action. Let $K_1,K_2$ be Killing fields generating the $T^2$-action, and let $\mathcal{G}=(\mathcal{G}_{ij})=g(K_{i},K_{j})$ be the Gram matrix. 

A number of different regions in $M$ are distinguished according to the possible stabilisers of the $T^2$-action. In regions where ${\rm rank}(\mathcal{G})=2$ the action is free and points have trivial stabiliser, whereas in regions with ${\rm rank}(\mathcal{G})<2$ the action degenerates, and $\det \mathcal{G} =0$. If ${\rm rank}(\mathcal{G})=1$, the action has an $S^1$-stabiliser; the corresponding region is called a {\it rod} and represents a $T^2$-invariant divisor in $M$. Points with ${\rm rank}(\mathcal{G})=0$ have stabiliser the whole $T^2$, and are called {\it turning points} or fixed points. 

The regularity of $(M,g)$ is analysed by imposing certain conditions on the above data. If the manifold is non-compact and the $T^2$-action has $n$ fixed points $p_{1},...,p_{n}$, one has a collection of $n+1$ rods $R_{k}$, $k=0,...,n$, of which $n-1$ are `finite' rods (non-contractible 2-spheres in $M$) and $2$ are `semi-infinite' rods (punctured 2-spheres). Two consecutive rods $R_{k}$, $R_{k+1}$ intersect at the fixed point $p_k$. Associated to every rod $R_{k}$ there is a generator $v_{k}$ of $T^2$ that vanishes on $R_{k}$. Suppose that $e_1,e_2$ are generators of $T^2$ which have closed orbits with periodicity $2\pi$. Writing $v_{k}=a_{k}e_1+b_{k}e_2$, the orbits of $v_{k}$ will be closed and have periodicity $2\pi$ iff $a_k$ and $b_k$ are coprime integers. If this is the case, then there will be no conical singularity along the rod $R_k$. However, there might still be an orbifold singularity at the intersection point $p_{k}$ of rods $R_k$ and $R_{k+1}$. The extra condition needed for the absence of an orbifold singularity at $p_{k}$ is $a_{k}b_{k+1}-b_{k}a_{k+1} = \pm 1$, where recall $(a_k,b_k)$ and $(a_{k+1},b_{k+1})$ are already assumed to be coprime integers. If the above conditions are satisfied for all rods and turning points, then there will be no conical or orbifold singularities \cite{CT0}.

Equivalently, we can analyse the absence of conical and orbifold singularities as follows. We first fix the scale of $v_{k}$ by requiring
\begin{align}\label{conical}
\lim_{{\rm Rod}\, R_k}\frac{|\d(|v_{k}|^{2})|^{2}}{4|v_{k}|^2}=1,
\end{align}
where the norms are computed with respect to the metric $g$, and the limit is taken in 
appropriate coordinates adapted to the rods. The Killing vector $v_{k}$ that vanishes on $R_{k}$ and satisfies \eqref{conical} is called the rod vector of $R_{k}$, and is unique up to sign. (If $v_k$ has closed orbits, then \eqref{conical} ensures that they are $2\pi$-periodic.) We then impose the following condition for the rod vectors of three consecutive rods:
\begin{align}\label{compatibility}
v_{k-1}+\varepsilon_{k} v_{k+1} = \ell_{k} v_{k},
\end{align}
where $\varepsilon_{k}=\pm1$ and $\ell_{k}\in\mathbb{Z}$. The absence of conical and orbifold singularities is then equivalent to the fulfilling of the conditions \eqref{conical} and \eqref{compatibility} for all rods.

In the case of ASD toric Einstein spaces, we saw in \S\ref{sec:multipoles} that the metric can be described in Tod form by \eqref{Tod1}-\eqref{Tod2}. The coordinates $(\rho,\zeta)$ are especially adapted to the rod structure: the turning points correspond to $(\rho,\zeta)=(0,\zeta_k)$, and the rods are $R_{k}=\{(\rho,\zeta)\,|\,\rho=0, \,\zeta_k<\zeta<\zeta_{k+1}\}$. The (normalised) rod vectors have been computed in \cite{ALR26}:
\begin{align}\label{rodvectors}
v_{k} = \begin{cases}
f'_{k}(-\omega_{k}\partial_{\psi}+\partial_{Y}) & \quad \text{if} \quad f'_{k}\neq 0 \\
\ell^2 f_{k}\partial_{\psi} & \quad \text{if} \quad f'_{k} = 0
\end{cases}
\end{align}
where $\omega_{k}:=\omega|_{R_{k}}$, $f_{k}:=f(\zeta_{k})$, $f'_{k}:=f'(\zeta)|_{R_{k}}$ are constants, and $f(\zeta)$ is a piecewise linear function defined by
\begin{align}\label{deff}
f(\zeta):=\hat{Z}|_{\rho=0}.
\end{align}
The facts that $f(\zeta)$ is piecewise linear and $\omega_{k}$ is constant have been shown in \cite{ALR26}. The regularity analysis then reduces to imposing the condition \eqref{compatibility} with the rod vectors given by \eqref{rodvectors}.

\subsection{Rod structure}
\label{sec:rods}
Note that for asymptotically hyperbolic manifolds, since there is a conformal boundary at $\hat{Z}=0$, some of the turning points and rods may actually be beyond the boundary. Therefore, the number of regularity conditions \eqref{compatibility} to be imposed depends on the locations of the boundary and the fixed points. 

For the solution \eqref{ASDCTintro}, taking the Killing fields $K_1=\partial_{\tau}$, $K_2=\partial_{\phi}$, the Gram matrix and its determinant are
\begin{align}
\mathcal{G} = 
\begin{pmatrix} 
\frac{F}{(x+\nu y)H} & \frac{G}{(x+\nu y)H} \\ 
\frac{G}{(x+\nu y)H} & \frac{G^2}{(x+\nu y)HF} - \frac{H\X\Y}{(x+\nu y)^3F}
\end{pmatrix}, 
\qquad 
\det\mathcal{G} = \frac{-\X\Y}{(x+\nu y)^4}.
\end{align}
Recalling the explicit form of $F,G,H$, we see that $(i)$ ${\rm rank}(\mathcal{G})=1$ when $\X=0$ {\it or} $\Y=0$, and $(ii)$ ${\rm rank}(\mathcal{G})=0$ when $\X=0$ {\it and} $\Y=0$; thus the roots of $\X$ (or $\Y$) characterise the degeneration of the $T^2$-action. With respect to the domain $\{r_2<x<r_3\}\cap\{r_1<y<r_2\}$ specified in \S\ref{sec:global}, we then see that there are in principle four rods: $\tilde{R}_1=\{x=r_2, \, r_1<y<r_2\}$, $\tilde{R}_2=\{r_2<x<r_3, \, y=r_1\}$, $\tilde{R}_3=\{x=r_3, \, r_1<y<r_2\}$, $\tilde{R}_4=\{r_2<x<r_3, \, y=r_2\}$. Furthermore, ${\rm rank}(\mathcal{G})=0$ at precisely four points: $(r_2,r_2), (r_2,r_1), (r_3,r_1), (r_3,r_2)$. We illustrate this in the diagram below.

\medskip
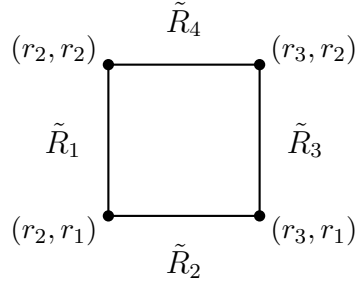
\begin{center}
\begin{tikzpicture}[scale=2,
    dot/.style={circle, fill=black, inner sep=1.5pt},
    dashedline/.style={dashed}
]

\draw[thick] (0,0) -- (1,0) -- (1,1) -- (0,1) -- (0,0);
\node[dot] at (0,0) {};
\node[dot] at (1,0) {};
\node[dot] at (1,1) {};
\node[dot] at (0,1) {};

\node at (-0.3,0.5) {$\tilde{R}_1$}; 
\node at (0.5,-0.3) {$\tilde{R}_2$}; 
\node at (1.3,0.5) {$\tilde{R}_3$}; 
\node at (0.5,1.3) {$\tilde{R}_4$}; 

\node at (-0.36,-0.1) {\small $(r_2,r_1)$}; 
\node at (1.36,-0.1) {\small $(r_3,r_1)$}; 
\node at (-0.36,1.1) {\small $(r_2,r_2)$}; 
\node at (1.36,1.1) {\small $(r_3,r_2)$}; 

\end{tikzpicture}
\end{center}

However, since the conformal infinity line $x+\nu y=0$ cuts through the coordinate domain, and each of the two Einstein metrics is defined on only one side of the boundary, the solutions will have $n$ and $4-n$ turning points for $n=1,2,3$ (depending on the value of $\nu$). Correspondingly, the rod $R_{k}$ will either coincide with $\tilde{R}_{k}$ or be a subset $R_{k}\subset\tilde{R}_{k}$, depending on whether $x+\nu y=0$ cuts through $\tilde{R}_{k}$ or not; see Fig. \ref{fig:solutionR} below. Nevertheless, the rod vectors for $\tilde{R}_{k}$ and $R_{k}$ are the same.

Note also that, since conformal infinity passes through exactly two rods, and by definition some circle collapses at a rod, the $T^2$-action on the boundary has exactly two singular orbits, which is what we claimed in \S\ref{sec:boundary}.

To compute the rod vectors \eqref{rodvectors}, we use the results of \S\ref{sec:multipoles}. The function $f(\zeta)$ in \eqref{deff} can be deduced from \eqref{4pole}:
\begin{align}\label{f4pole}
f(\zeta)= A+\alpha_1|\zeta-\zeta_1|+\alpha_2|\zeta-\zeta_2|+\alpha_3|\zeta-\zeta_3|,
\end{align}
where $A,\alpha_1,...,\alpha_3,\zeta_1,...,\zeta_3$ are given in \eqref{coeff4pole}-\eqref{coeff4pole2} in terms of $r_1,...,r_4,a_4,\nu$. Using \eqref{f4pole} we can compute the slopes $f'_1,...,f'_4$: 
\begin{equation}\label{slopes}
\begin{aligned}
f'_{1} ={}& -\alpha_1-\alpha_2-\alpha_3 =  \frac{-2r_2\nu(1+\nu)}{a_4 r_{12}r_{23}r_{24}}, \\
f'_{2} ={}& +\alpha_1-\alpha_2-\alpha_3 =  \frac{-2r_1\nu(1+\nu)}{a_4 r_{12}r_{13}r_{14}}, \\
f'_{3} ={}& +\alpha_1+\alpha_2-\alpha_3 = \frac{2r_3\nu(1+\nu)}{a_4 r_{13}r_{23}r_{24}}, \\
f'_{4} ={}& +\alpha_1+\alpha_2+\alpha_3 =  -f'_1,
\end{aligned}
\end{equation}
where $r_{ij}=r_{i}-r_{j}$. Note that we can choose $a_4$ to set $f'_1=1$ and $f'_4=-1$ as in \cite{ALR26}; we will not need to do that at the moment (cf. however \eqref{speciala4} below).

From \eqref{rodvectors}, we also need the constants $\omega_{k}$. This can be computed using \eqref{CTvariables}, which gives
\begin{align}\label{omegak}
\omega_{k} = \frac{k}{(\nu+1)^2}\left.\frac{G}{F}\right|_{R_{k}}. 
\end{align}

\subsection{Instantons}
\label{sec:regularity}

For concreteness, let us assume the coordinate domain illustrated in Fig.~\ref{FigureSDECT}. This leads to the diagrams in Fig. \ref{fig:diagrams}. We will focus on the regularity analysis of the Einstein metric labeled `R' in those diagrams, as this is globally conformal to the ALE K\"ahler metric from \S\ref{sec:kahler}. The possible diagrams are illustrated in Fig. \ref{fig:solutionR}.

\medskip
\begin{center}
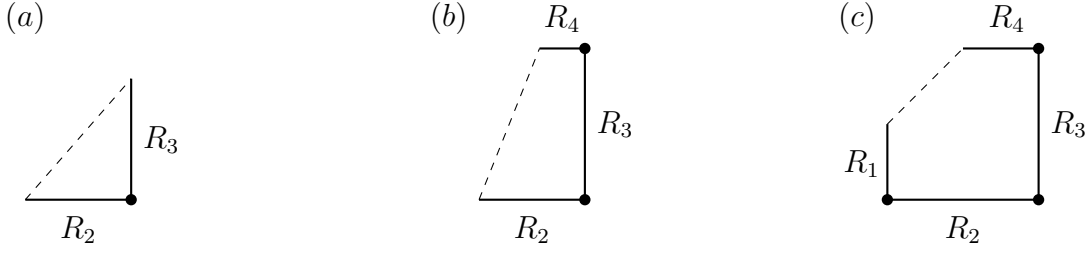
\begin{figure}
\begin{tikzpicture}[scale=2,
    dot/.style={circle, fill=black, inner sep=1.5pt},
    dashedline/.style={dashed}
]

\begin{scope}[xshift=0cm]
\draw[thick] (0.3,0) -- (1,0) -- (1,0.8);
\draw[dashedline] (0.3,0) -- (1,0.8);
\node[dot] at (1,0) {};
\node at (0.3,1.2) {$(a)$};
\node at (0.65,-0.2) {$R_2$};
\node at (1.2,0.4) {$R_3$};
\end{scope}

\begin{scope}[xshift=3cm]
\draw[thick] (0.3,0) -- (1,0) -- (1,1) -- (0.7,1);
\draw[dashedline] (0.3,0) -- (0.7,1);
\node[dot] at (1,0) {};
\node[dot] at (1,1) {};
\node at (0.1,1.2) {$(b)$};
\node at (0.65,-0.2) {$R_2$};
\node at (1.2,0.5) {$R_3$};
\node at (0.85,1.2) {$R_4$};
\end{scope}

\begin{scope}[xshift=6cm]
\draw[thick] (0,0.5) -- (0,0) -- (1,0) -- (1,1) -- (0.5,1);
\draw[dashedline] (0,0.5) -- (0.5,1);
\node[dot] at (0,0) {};
\node[dot] at (1,0) {};
\node[dot] at (1,1) {};
\node at (-0.2,1.2) {$(c)$};
\node at (-0.17,0.24) {$R_1$};
\node at (0.5,-0.2) {$R_2$};
\node at (1.2,0.5) {$R_3$};
\node at (0.79,1.2) {$R_4$};
\end{scope}

\end{tikzpicture}
\caption{Diagrams corresponding to the solution `R' in Fig. \ref{fig:diagrams}.}
\label{fig:solutionR}
\end{figure}
\end{center}

Consider first case $(a)$. We see that the geometry has only one fixed point and two semi-infinite rods, $R_2$ and $R_3$. Thus in this case the extra regularity condition \eqref{compatibility} is trivial. In order to avoid conical singularities along the rods, the orbits of the rod vectors $v_2,v_3$ should be independently identified with period $2\pi$, that is
\begin{align}
(\psi,Y)\sim(\psi - 2\pi f'_2\omega_2,Y+2\pi f'_2), \qquad
(\psi,Y)\sim(\psi - 2\pi f'_3\omega_3,Y+2\pi f'_3).
\end{align}
With these identifications, the solution is regular. We did not need to fix any of the original effective parameters, so this gives a 4-parameter gravitational instanton. Note that the topology is that of the 4-ball, $B^4$.

Consider now case $(b)$. The solution has two fixed points, two semi-infinite rods, and one finite rod. The regularity requirement \eqref{compatibility} reduces to one condition: $v_2+\varepsilon v_4 = p v_3$, with $\varepsilon=\pm1$, $p\in\mathbb{Z}$. Using the rod vectors \eqref{rodvectors} (where assume $f'_k\neq0$ for concreteness), the conditions are 
\begin{equation}\label{regcaseb}
\begin{aligned}
f'_2+\varepsilon f'_4 - p f'_3 ={}&0 , \\
f'_2\omega_2+\varepsilon f'_4\omega_4 - p f'_3\omega_3 ={}&0.
\end{aligned}
\end{equation}
This can be expressed explicitly in terms of the parameters $r_1,...,r_4,a_4,\nu$ using \eqref{slopes} and \eqref{omegak}. One can then check that a possible solution to the conditions \eqref{regcaseb} is $\varepsilon=+1$ and
\begin{equation}
\begin{aligned}
r_1={}& -\frac{r_3 (-1 + p + \nu) (r_2 + (-2 + p) r_3 + (-1 + p) r_2 \nu)}{(1 + (-1 + p) \nu) ((-2 + p) r_2 \nu + r_3 (-1 + p + \nu))}, \\ 
r_4 ={}& -\frac{r_3 (-1 + p + \nu)}{1 + (-1 + p) \nu}.
\end{aligned}
\end{equation}
The parameters are now $r_2,r_3,a_4,\nu$ and an integer $p$. But since $r_2$ and $a_4$ are pure gauge, the solution depends effectively on three parameters, of which two are continuous and one is an integer. Topologically, the finite rod gives a non-contractible 2-cycle, and the manifold is the total space of the line bundle $\mathcal{O}(-p)\to\mathbb{CP}^1$, or equivalently, the toric resolution of the singularity $\mathbb{C}^2/\mathbb{Z}_p$ where $\mathbb{Z}_p$ acts by scalar multiplication \cite{CS}. The boundary is the lens space $L(p,1)$. One can also check that the case $p=0$ does not solve the regularity conditions \eqref{regcaseb}, so we must have $p\neq0$.

Finally, consider case $(c)$. We see that the solution has three fixed points, two semi-infinite rods, and two finite rods. The regularity requirement \eqref{compatibility} reduces to two conditions: $v_1+\varepsilon_2 v_3 = \ell_2 v_2$ and 
$v_2+\varepsilon_3 v_4 = \ell_3 v_3$, with $\varepsilon_2,\varepsilon_3=\pm1$ and $\ell_2, \ell_3 \in\mathbb{Z}$. In components:
\begin{subequations}\label{regcasec}
\begin{align}
f'_1+\varepsilon_2 f'_3 - \ell_2 f'_2 ={}&0, \label{cond1} \\
f'_1\omega_1+\varepsilon_2 f'_3\omega_3 - \ell_2 f'_2\omega_2={}&0, \label{cond2} \\ 
f'_2+\varepsilon_3 f'_4 - \ell_3 f'_3 ={}&0, \label{cond3} \\
f'_2\omega_2+\varepsilon_3 f'_4\omega_4 - \ell_3 f'_3\omega_3={}&0. \label{cond4}
\end{align}
\end{subequations}
To analyse these conditions, we express them in terms of $r_1,...,r_4,a_4,\nu$ using \eqref{slopes} and \eqref{omegak}, and try to see if the parameters can be chosen so as to satisfy \eqref{regcasec}. 
For example, one can check that \eqref{cond1} and \eqref{cond2} are solved for $\varepsilon_2=+1$ and
\begin{equation}\label{r3r4}
\begin{aligned}
r_3 ={}& -\frac{r_1 (1 + (-1 + \ell_2) \nu) ((-2 + \ell_2) r_1 \nu + r_2 (-1 + \ell_2 + \nu))}{(-1 + \ell_2 + \nu) (r_1 + (-2 + \ell_2) r_2 + (-1 + \ell_2) r_1 \nu)}, \\
r_4 ={}& \frac{r_1 (-1 + \nu - \ell_2 \nu)}{-1 + \ell_2 + \nu}.
\end{aligned}
\end{equation}
Replacing these values of $r_3,r_4$ into \eqref{cond3}, we set $\varepsilon_3=+1$ and choose $r_1$ so as to satisfy the equation. The expression of $r_1$ in terms of $\nu,\ell_2,\ell_3$ that we get is long and complicated, so we will not display it here. We replace this back into \eqref{r3r4}, so we end up with $r_1=r_1(\nu,\ell_2,\ell_3)$, $r_3=r_3(\nu,\ell_2,\ell_3)$ and $r_4=r_4(\nu,\ell_2,\ell_3)$. The next step is then to replace these functions of $(\nu,\ell_2,\ell_3)$ into \eqref{cond4}. This gives a notoriously long and cumbersome rational function of $(\nu,\ell_2,\ell_3)$ involving many roots and different powers of $\nu,\ell_2,\ell_3$, so we will not display it here either. The last step is to solve the resulting condition for one of the remaining non-trivial parameters, which are $\nu,\ell_2,\ell_3$. Since $\ell_2,\ell_3$ are required to be integers, it is convenient for solve for $\nu$. While solutions exist, we have not been able to get a closed-form expression for $\nu$ as a function $\nu(\ell_2,\ell_3)$ in the most general case. 

The end result is a list of functions $\nu(\ell_2,\ell_3)$ and $r_i(\ell_2,\ell_3)$ for $i=1,3,4$ that solve the regularity conditions \eqref{regcasec} (recall that $r_2$ and $a_4$ are pure gauge). However, one should still check that all consistency conditions are indeed satisfied. Namely, the diagram $(c)$ assumes $-1<\nu<0$ and that three of the turning points are inside the domain $\hat{Z}>0$. In particular, the solution should satisfy the conditions listed in \cite[Prop. 4.6]{ALR26}, that is $f_k>0$, $\alpha_k<0$ \footnote{The parameters $\alpha_k$ here correspond to $a_k$ in \cite{ALR26}.}, and $A>0$ (where we note that this assumes that the slopes at infinity are normalised as $f'_1=+1$, $f'_4=-1$). Given the complexity of the functions $\nu(\ell_2,\ell_3)$, $r_i(\ell_2,\ell_3)$, we will not attempt a general analysis, but simply show explicitly that solutions exist.

Indeed, consider the case 
\begin{align}
\ell_3=\ell_2 \equiv p,
\end{align}
where $p$ will be required to be an integer. Following the steps described above, the equation \eqref{cond4} simplifies significantly and becomes
\begin{equation}
\begin{aligned}
\frac{3\nu}{2p\lambda} & \left[ 2 + p - p^2 + (-2 + p) (1 + p) \nu \right. \\
 & \left. + \sqrt{(-1 + p) (-4 (-1 + \nu)^2 + 3 p^2 (-1 + \nu)^2 + p^3 (1 + \nu (6 + \nu)))} \right] = 0.
\end{aligned}
\end{equation}
Since in our current situation $\nu=0$ is not allowed, the term inside the brackets must vanish. This is quadratic in $\nu$ so there are two solutions. We choose the following solution:
\begin{align}
\nu = \frac{-2 + 2 p^2 - p^3 + \sqrt{(-2 + p)^2 p^2 (-1 + p^2)}}{-2 + p^2}
\end{align}
For $p>\sqrt{2}$, this function is continuous, negative-definite, monotonically increasing, and it goes from $\nu=-1$ at $p=2$ to $\nu=0$ for $p\to\infty$. Since the value $\nu=-1$ is not allowed, we must take $p\geq3$. So we indeed have $-1<\nu<0$ and $p\in\mathbb{Z}$, and the regularity conditions are satisfied.

As an example, take $p = 3$. Then, one can check that \eqref{regcasec} is satisfied for
\begin{align}
\nu = \frac{1}{7}(-11+6\sqrt{2}), \qquad r_1=(3+2\sqrt{2})r_2, \qquad r_{3}=(3-2\sqrt{2})r_2, \qquad r_4=-r_2
\end{align}
(where one can choose e.g. $r_2=-1$, which gives $\nu\approx-0.36$, $r_1\approx-5.82$, $r_2=-1$, $r_3\approx-0.17$ and $r_4=1$).
The slopes at infinity, $f'_1, f'_4$, can be set to $f'_1=1=-f'_4$ by choosing 
\begin{align}\label{speciala4}
a_4=\frac{-58+45\sqrt{2}}{98 r_2^2}.
\end{align}
This gives $\alpha_1 = \frac{-3}{8}$, $\alpha_2 = \frac{-1}{4}$, $\alpha_3 = \alpha_1$, $f_1=\frac{(3+\sqrt{2})}{7}$, $f_2 = \frac{3(3+\sqrt{2})}{14}$, $f_3 = f_1$, and  $A = \frac{3(3+\sqrt{2})}{7}$. All regularity and consistency conditions are then satisfied, so we indeed a gravitational instanton.

For other integer values of $\ell_2,\ell_3$ (with $\ell_2>2,\ell_3>2$, cf. \cite{ALR26}), the analysis is similar. Case $(c)$ then gives a 2-parameter gravitational instanton, where the parameters are the integers $\ell_2,\ell_3$. The two finite rods give two non-contractible 2-cycles, and the topology of $M$ is the same as that in \cite{CS} for $b_2(M)=2$ (i.e, a toric resolution of $\mathbb{C}^2/\Gamma$ where the exceptional divisor has two 2-spheres). The integers $\ell_2,\ell_3$ give the self-intersection numbers of the two 2-cycles.

This finishes the proof of Theorem \ref{thm:instantons}.

\section{Final remarks}
\label{sec:discussion}

We conclude with some comments about the generality of the ansatz in \S\ref{sec:ansatz}, especially for the case $\lambda>0$. As mentioned in \S\ref{sec:introduction}, it is well-known that toric K\"ahler structures can be encoded in Delzant polytopes, cf. \cite{Abreu}. For Hermitian-Einstein instantons with $\lambda>0$, LeBrun's classification \cite{LeBrun2} shows that the possible polytopes are the following:
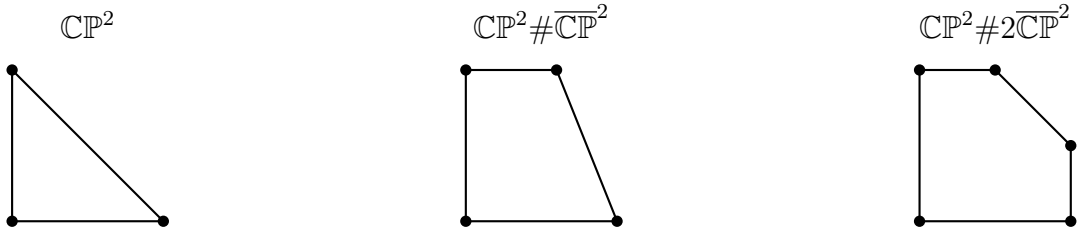
\begin{center}
\begin{tikzpicture}[scale=2]
\begin{scope}[xshift=0cm]
\node[draw=none, fill=none] at (0.5,1.3) {$\mathbb{CP}^2$};
\coordinate (A) at (0,0);
\coordinate (B) at (0,1);
\coordinate (C) at (1,0);
\draw[thick] (A) -- (B) -- (C) -- cycle;
\node[circle, fill=black, inner sep=1.5pt] at (A) {};
\node[circle, fill=black, inner sep=1.5pt] at (B) {};
\node[circle, fill=black, inner sep=1.5pt] at (C) {};
\end{scope}
\begin{scope}[xshift=3cm]
\node[draw=none, fill=none] at (0.5,1.3) {$\mathbb{CP}^2 \# \overline{\mathbb{CP}}^2$};
\coordinate (A) at (0,0);
\coordinate (B) at (0,1);
\coordinate (C) at (0.6,1);
\coordinate (D) at (1,0);
\draw[thick] (A) -- (B) -- (C) -- (D) -- cycle;
\node[circle, fill=black, inner sep=1.5pt] at (A) {};
\node[circle, fill=black, inner sep=1.5pt] at (B) {};
\node[circle, fill=black, inner sep=1.5pt] at (C) {};
\node[circle, fill=black, inner sep=1.5pt] at (D) {};
\end{scope}
\begin{scope}[xshift=6cm]
\node[draw=none, fill=none] at (0.5,1.3) {$\mathbb{CP}^2 \# 2\overline{\mathbb{CP}}^2$};
\coordinate (A) at (0,0);
\coordinate (B) at (0,1);
\coordinate (C) at (0.5,1);
\coordinate (D) at (1,0.5);
\coordinate (E) at (1,0);
\draw[thick] (A) -- (B) -- (C) -- (D) -- (E) -- cycle;
\node[circle, fill=black, inner sep=1.5pt] at (A) {};
\node[circle, fill=black, inner sep=1.5pt] at (B) {};
\node[circle, fill=black, inner sep=1.5pt] at (C) {};
\node[circle, fill=black, inner sep=1.5pt] at (D) {};
\node[circle, fill=black, inner sep=1.5pt] at (E) {};
\end{scope}
\end{tikzpicture}
\end{center}
\medskip
The polytopes represent the image of the moment map of the $T^2$-action, and can be identified with the orbit space described in moment-map coordinates $x_1,x_2$. The vertices have $T^2$-stabiliser, and the edges have $S^1$-stabiliser, or in the language of \S\ref{sec:rodstructure}, the vertices are turning points and the edges are rods. This means that the determinant of the Gram matrix of toric Killing fields must vanish on the edges. 

In terms of the Toda variables of \S\ref{sec:THE}, we see from \eqref{Todavariables} that this determinant is precisely the function $e^{u}$ that is required to solve the Toda equation \eqref{MTE}. Now recall that for the ansatz in \S\ref{sec:ansatz}, using \eqref{eu}, \eqref{ABansatz}, \eqref{CDansatz} we have 
\begin{align}\label{eudiscussion}
e^{u} = -\frac{\X\Y}{(x-y)^4}
\end{align}
where the functions $\X=\X(x)$ and $\Y=\Y(y)$ are yet to be determined by solving the Toda equation. Although the coordinates $(x,y)$ are generically not moment-map coordinates $(x_1,x_2)$ (in particular, the former are orthogonal while the latter are not), both sets of coordinates capture the degeneration of the $T^2$-action, which is an invariant. 

Now, the fact that $\X$ depends only on $x$ and $\Y$ depends only on $y$ implies that the function \eqref{eudiscussion} can only vanish on vertical and horizontal lines in the $(x,y)$-plane. From this we deduce, without solving the Einstein equations, that \eqref{eudiscussion} can only describe the first two instantons above. That is, the $(x,y)$-coordinates can be chosen such that: for $\CP^2$, one of the vertices is at infinity, and for $\CP^2\#\overline{\CP}{}^{2}$, the trapezoid is transformed to a rectangular shape. Thus both $\CP^2$ and $\CP^2\#\overline{\CP}{}^{2}$ can be described with \eqref{eudiscussion}; in particular, they should be obtained as special limits of the Pleba\'nski-Demia\'nski family with $\lambda>0$.

However, we now also see that the pentagon corresponding to $\CP^2\#2\overline{\CP}{}^{2}$ cannot be described by the ansatz \eqref{eudiscussion}, since it cannot be transformed to a rectangular shape without removing a vertex. This provides an explanation for why \eqref{CTansatz}-\eqref{CTansatzfunctions} does not lead to the CLW metric, without needing to solve the Einstein equations\footnote{In the language of rod structures, the ansatz \eqref{eudiscussion} can describe up to four turning points and four rods, whereas the CLW metric has five turning points and five rods.}.

\smallskip
It is also worth briefly commenting on the similarity between the above polytopes and Fig. \ref{fig:solutionR}. Namely, both sets of diagrams codify the degeneration of the torus action (in different orbit-space coordinates), and the diagrams in Fig. \ref{fig:solutionR} look (after a projective transformation) like the polytopes with one of the edges removed. This is reminiscent of what happens in the Ricci-flat case, where the (Hermitian) ALF instantons with $\lambda=0$ correspond to the compact instantons with $\lambda>0$ after removing an edge of the polytope \cite{BG}. However, the difference here is that, in the Ricci-flat case, the removed edge corresponds to a 2-sphere fixed by the $T^2$-action (i.e. a rod), which results naturally in ALF asymptotics, cf. \cite{OSD}. In fact, as shown by Li \cite{Li23}, the ALF Ricci-flat instantons can be compactified by adding a 2-sphere at infinity.
By contrast, in Fig. \ref{fig:solutionR} the ``removed edge'' corresponds to the 3-dimensional conformal boundary (which is not a rod), and the conformal compactification produces an ALE metric, which is compactified by adding an (orbifold) point at infinity instead of a sphere. Note that this difference is also manifested in the fact that the topology of the instantons in \S\ref{sec:regularity} is different from the topology of the ALF Ricci-flat instantons.

\smallskip
Finally, we expect the ASD Einstein metric studied in this work to correspond to an ASD limit of a non-(anti)self-dual Hermitian-Einstein metric with $\lambda<0$. In contrast to the case $\lambda>0$, much less is known about the non-ASD, $\lambda<0$ case, cf. \cite{Li:2025yll}. The solution should be obtained by relaxing the assumption \eqref{AnsatzX}, which, from the above discussion, should produce a metric with more than four rods. We leave this for future work.

\subsection*{Acknowledgements}
I am very grateful to Tim Adamo, Lars Andersson, Yu Chen, Maciej Dunajski, Claude LeBrun, Mingyang Li, James Lucietti, Rosa Sena-Dias, and Paul Tod for many enlightening conversations. I also gratefully acknowledge the hospitality and support of the Simons Center for Geometry and Physics, Stony Brook University, New York, during the program `Einstein 4-manifolds and Gravitational Instantons' in the Winter 2026, at which a major part of the research for this paper was performed. The author is supported via the ERC Consolidator/UKRI Frontier grant TwistorQFT EP/Z000157/1.

\end{document}